\documentclass[11pt]{article}
\RequirePackage{amsthm,amsmath, amssymb}
 \usepackage{epsfig}
 \usepackage{graphicx}
 \usepackage{pict2e}
\newcommand{\bx}{\mathbf{x}}
\newcommand{\bX}{\mathbf{X}}

\newcommand{\bF}{\mathbf{F}}
\newcommand{\CB}{\mathbf{CB}}
\newcommand{\bdf}{\mathbf{f}}
\newcommand{\Reals}{\mathbb{R}}
\newcommand{\Ints}{\mathbb{Z}}
\newcommand{\Ex}{\mathbb{E}}
\renewcommand{\Pr}{\mathbb{P}}
 \newcommand{\sfrac}[2]{{\textstyle\frac{#1}{#2}}}
\newcommand{\FF}{\mathcal{F}}
\newcommand{\GG}{\mathcal{G}}
\newcommand{\HH}{\mathcal{H}}
\newcommand{\II}{\mathcal{I}}

\newcommand{\cd}{\stackrel{d}{\to}}
\def\ind{{\rm 1\hspace{-0.90ex}1}}
\newcommand{\Fspace}{\mathbb{F}}

\newcommand{\init}{\mathrm{init}}

\newcommand{\rmod}{\mathrm{modi}}
\newtheorem{Lemma}{Lemma}
 \newtheorem{Proposition}[Lemma]{Proposition}
 
 \newtheorem{Theorem}[Lemma]{Theorem}
 \newtheorem{Conjecture}[Lemma]{Conjecture}

 \newtheorem{OP}[Lemma]{Open Problem}

 \newtheorem{Corollary}[Lemma]{Corollary}

 \title{Waves in a Spatial Queue: Stop-and-Go at Airport Security}
 \author{David J. Aldous\thanks{Research supported by NSF Grants DMS-106998 and 1504802.}\\
Department of Statistics\\
367 Evans Hall \#\  3860\\
U.C. Berkeley CA 94720\\
aldous@stat.berkeley.edu}

\begin{document}
\maketitle

\begin{abstract}
We model a long queue of humans by a continuous-space model in which, when a customer moves forward, 
they stop a random distance behind the previous customer, but do not move at all if their distance behind the previous 
customer is below a threshold.  The latter assumption leads to ``waves" of motion in which
only some random number $W$ of customers move.  
We prove that $\Pr(W > k)$ decreases as order $k^{-1/2}$; in other words, for large $k$ the $k$'th customer moves 
on average only once every order $k^{1/2}$ service times. 
A more refined analysis relies on a non-obvious asymptotic relation to the 
coalescing Brownian motion process; we give a careful outline of such an analysis without attending to all the technical details.
\end{abstract}



\section{Introduction}
Imagine you are the 100th person in line at an airport security checkpoint. 
As people reach the front of the line they are  being processed steadily, at rate 1 per unit time.
But you move less frequently, and when you do move, you typically move several units of distance, 
where 1 unit distance is the average distance between successive people standing  in the line.

This phenomenon is easy to understand qualitatively. 
When a person leaves the checkpoint, the next person moves up to the checkpoint, the next person  moves up and stops behind 
the now-first person, and so on, but this ``wave" of motion often does not extend through the entire long line; instead, 
some person will move only a short distance, and the person behind will decide not to move at all.
Intuitively, when you are around the $k$'th position in line,  there must be some number $a(k)$ representing both the
 average time between your moves and the average distance you move when you do move -- these are equal because
you are moving forwards at average speed $1$.  
In other words, the number $W$ of people who move at a typical step has distribution 
$\Pr(W \ge k) = 1/a(k)$.
This immediately suggests the question of  how fast  $a(k)$ grows with $k$.
In this paper we will present a stochastic model and prove that, in the model,  $a(k)$ grows as order $k^{1/2}$.

In classical {\em queueing theory} \cite{kleinrock1975queueing},
 randomness enters via assumed randomness 
of arrival and service times.  In contrast, 
 even though we are modeling a literal queue, randomness in our model arises in a quite different way, via
 each customer's choice of exactly how far behind the preceding  customer  they choose to stand, after each move.  
That is, we assume that ``how far behind" is chosen 
(independently for each person and time) 
from a given 
probability measure $\mu$ on an interval $[c_-,c^+]$ where $c_- > 0$. 
We interpret this interval as a ``comfort zone" for proximity to other people. 
By scaling we may assume $\mu$ has mean $1$, and then (excluding the deterministic case) 
$\mu$ has some variance $0 < \sigma^2 < \infty$.

In words, the model is 
\begin{quote}
when the person in front of you moves forward to a new position, then you move to a new position 
at a random distance (chosen from distribution  $\mu$) behind them, unless their new position is less than distance $c^+$ in front of your existing position, in which case you don't move, and therefore nobody behind you moves.
\end{quote}
For the precise mathematical setup it is easiest to consider an infinite queue. 
A state $\bx = (x_0, x_1, x_2, \ldots)$ of the process is a configuration
\[ 0 = x_0 < x_1 < x_2 < x_3 \ldots, \quad x_i - x_{i-1} \in [c_-,c^+] \]
representing the positions of the customers on the half-line $[0,\infty)$.
The state space $\mathbb{X}$ is the set of such configurations.
The process we study -- let us call it the {\em spatial queue} process -- is the discrete-time 
$\mathbb{X}$-valued Markov chain 
$(\bX(t), t = 0,1,2,\ldots)$  whose transitions from an arbitrary initial state $\bx$ to $\bX(1)$ have the following distribution.
Take i.i.d. $(\xi_i, i \ge 1)$ with distribution $\mu$ and define
\begin{eqnarray*}
X_0(1) &=& 0 \\
W(1) &=& \min \{i: \xi_1 + \ldots + \xi_i \ge  x_{i+1} - c^+\} \\
X_i(1) &=&   \xi_1 + \ldots + \xi_i , \ 1 \le  i < W(1)\\
&=& x_{i+1}, \ i \ge W(1).
\end{eqnarray*}
This simply formalizes the verbal description above.
Note $W(1)$ is the number of customers who moved (including the customer moving to the service position $0$) 
in that step; call $W(1)$ the {\em wave length} at step $1$.
As more terminology: in configuration $\bx = (x_i, i \ge 0)$ the customer at {\em position } $x_i$ has 
{\em rank} $i$.
So in each step, a customer's rank will decrease by $1$ whereas their position may stay the same or may decrease by some 
distance in $(0,\infty)$; recall that the random variables $\xi$ do not directly indicate the distance a customer moves, 
but instead indicate the distance at which the customer stops relative to the customer ahead.
For a general step of the process we take 
i.i.d. $(\xi_i(t), i \ge 1, t \ge 1)$ with distribution $\mu$ and define
\begin{eqnarray}
W(t) &=& \min \{i: \xi_1(t) + \ldots + \xi_i (t) \ge  X_{i+1}(t-1) - c^+\} \label{Xit1}\\
X_i(t) &=&   \xi_1(t) + \ldots + \xi_i(t) , \  1 \le i < W(t) \label{Xit2}\\
&=& X_{i+1}(t-1), \ i \ge W(t). \label{Xit3}
\end{eqnarray}
Before continuing we should point out something we cannot prove, but strongly believe.
\begin{Conjecture}
\label{C-1}
Under appropriate ``non-lattice" conditions on $\mu$, 
the spatial queue process has a unique stationary distribution, $\bX(\infty)$ say, and from any initial 
configuration we have $\bX(t) \to_d \bX(\infty)$ as $t \to \infty$.
\end{Conjecture}
Understanding the space-time trajectory of a typical customer far away from the head of the queue 
is essentially equivalent to understanding the large values of the  induced process $(W(t), t \ge 1)$ of wave lengths, 
so that will be our focus.
Intuitively we study properties of the stationary distribution, or of 
a typical customer in the stationary process, but we are forced to phrase results in a time-averaged way.
Our fundamental result concerns the order of magnitude of long waves.
Define
\[ \rho^+(j) = \limsup_{\tau \to \infty} \tau^{-1} \sum_{t=1}^\tau \Pr(W(t) > j) , \ 
\rho_-(j) = \liminf_{\tau \to \infty} \tau^{-1} \sum_{t=1}^\tau \Pr(W(t) > j) . \]
\begin{Proposition}
\label{P-main}
There exist constants $0 < a_- < a^+ < \infty$ (depending on $\mu$) such that
\begin{equation}
 a_- j^{-1/2} \le \rho_-(j) \le \rho^+(j) \le a^+ j^{-1/2}, \ j \ge 1 . 
 \label{ajj}
 \end{equation}
 \end{Proposition}
In other words for a rank-$j$ customer, for large $j$, the average time between moves is order $j^{1/2}$. 
This result is proved in sections \ref{sec:UB} and \ref{sec:LB} as Proposition \ref{P:ub} and Corollary \ref{C:lb}.

Building upon this order-of-magnitude result (Proposition \ref{P-main}),
sharper results can be obtained by exploiting a  connection with the coalescing Brownian motion (CBM) process, 
which we start to explain next.

\subsection{Graphics of realizations}
\label{sec:graphics}
The graphics below not only  illustrate the qualitative behavior of the process but also, via a non-obvious 
representation, suggest  the main result.

\setlength{\unitlength}{0.05in}
\begin{picture}(100,67)
\put(0,0){\circle*{0.61}}
\put(12,0){\circle*{0.61}}
\put(22,0){\circle*{0.61}}
\put(29,0){\circle*{0.61}}
\put(36,0){\circle*{0.61}}
\put(50,0){\circle*{0.61}}
\put(63,0){\circle*{0.61}}
\put(70,0){\circle*{0.61}}
\put(79,0){\circle*{0.61}}
\put(86,0){\circle*{0.61}}
\put(94,0){\circle*{0.61}}
\put(22,0.5){\line(-15,4){15}}
\put(36,0.5){\line(0,1){4}}
\put(63,0.5){\line(0,1){4}}
\put(79,0.5){\line(0,1){4}}
\put(94,0.5){\line(0,1){4}}
\put(0,5){\circle*{0.61}}
\put(7,5){\circle*{0.61}}
\put(21,5){\circle*{0.61}}
\put(36,5){\circle*{0.61}}
\put(50,5){\circle*{0.61}}
\put(63,5){\circle*{0.61}}
\put(70,5){\circle*{0.61}}
\put(79,5){\circle*{0.61}}
\put(86,5){\circle*{0.61}}
\put(94,5){\circle*{0.61}}
\put(7,5.5){\line(-7,4){7}}
\put(36,5.5){\line(-16,4){16}}
\put(63,5.5){\line(-19,4){19}}
\put(79,5.5){\line(-7,4){7}}
\put(94,5.5){\line(0,1){4}}
\put(0,10){\circle*{0.61}}
\put(10,10){\circle*{0.61}}
\put(20,10){\circle*{0.61}}
\put(35,10){\circle*{0.61}}
\put(44,10){\circle*{0.61}}
\put(60,10){\circle*{0.61}}
\put(72,10){\circle*{0.61}}
\put(86,10){\circle*{0.61}}
\put(94,10){\circle*{0.61}}
\put(20,10.5){\line(-13,4){13}}
\put(44,10.5){\line(-19,4){19}}
\put(72,10.5){\line(-26,4){26}}
\put(94,10.5){\line(-33,4){33}}
\put(78,14.5){\line(22,-3){22}}
\put(0,15){\circle*{0.61}}
\put(7,15){\circle*{0.61}}
\put(17,15){\circle*{0.61}}
\put(25,15){\circle*{0.61}}
\put(34,15){\circle*{0.61}}
\put(46,15){\circle*{0.61}}
\put(55,15){\circle*{0.61}}
\put(61,15){\circle*{0.61}}
\put(69,15){\circle*{0.61}}
\put(78,15){\circle*{0.61}}
\put(91,15){\circle*{0.61}}
\put(100,15){\circle*{0.61}}
\put(7,15.5){\line(-7,4){7}}
\put(25,15.5){\line(0,1){4}}
\put(46,15.5){\line(0,1){4}}
\put(61,15.5){\line(0,1){4}}
\put(78,15.5){\line(0,1){4}}
\put(100,15.5){\line(0,1){4}}
\put(0,20){\circle*{0.61}}
\put(13,20){\circle*{0.61}}
\put(25,20){\circle*{0.61}}
\put(34,20){\circle*{0.61}}
\put(46,20){\circle*{0.61}}
\put(55,20){\circle*{0.61}}
\put(61,20){\circle*{0.61}}
\put(69,20){\circle*{0.61}}
\put(78,20){\circle*{0.61}}
\put(91,20){\circle*{0.61}}
\put(100,20){\circle*{0.61}}
\put(25,20.5){\line(-20,4){20}}
\put(46,20.5){\line(-26,4){26}}
\put(61,20.5){\line(-14,4){14}}
\put(78,20.5){\line(-17,4){17}}
\put(100,20.5){\line(-20,4){20}}
\put(97,24.5){\line(3,-1){3}}
\put(0,25){\circle*{0.61}}
\put(5,25){\circle*{0.61}}
\put(13,25){\circle*{0.61}}
\put(20,25){\circle*{0.61}}
\put(33,25){\circle*{0.61}}
\put(47,25){\circle*{0.61}}
\put(56,25){\circle*{0.61}}
\put(61,25){\circle*{0.61}}
\put(74,25){\circle*{0.61}}
\put(80,25){\circle*{0.61}}
\put(87,25){\circle*{0.61}}
\put(97,25){\circle*{0.61}}
\put(5,25.5){\line(-5,4){5}}
\put(20,25.5){\line(0,1){4}}
\put(47,25.5){\line(0,1){4}}
\put(61,25.5){\line(0,1){4}}
\put(80,25.5){\line(0,1){4}}
\put(97,25.5){\line(0,1){4}}
\put(0,30){\circle*{0.61}}
\put(10,30){\circle*{0.61}}
\put(20,30){\circle*{0.61}}
\put(33,30){\circle*{0.61}}
\put(47,30){\circle*{0.61}}
\put(56,30){\circle*{0.61}}
\put(61,30){\circle*{0.61}}
\put(74,30){\circle*{0.61}}
\put(80,30){\circle*{0.61}}
\put(87,30){\circle*{0.61}}
\put(97,30){\circle*{0.61}}
\put(20,30.5){\line(-14,4){14}}
\put(47,30.5){\line(-20,4){20}}
\put(61,30.5){\line(-17,4){17}}
\put(80,30.5){\line(-22,4){22}}
\put(97,30.5){\line(-21,4){21}}
\put(91,34.5){\line(8,-1){9}}
\put(0,35){\circle*{0.61}}
\put(6,35){\circle*{0.61}}
\put(17,35){\circle*{0.61}}
\put(27,35){\circle*{0.61}}
\put(35,35){\circle*{0.61}}
\put(44,35){\circle*{0.61}}
\put(50,35){\circle*{0.61}}
\put(58,35){\circle*{0.61}}
\put(69,35){\circle*{0.61}}
\put(76,35){\circle*{0.61}}
\put(84,35){\circle*{0.61}}
\put(91,35){\circle*{0.61}}
\put(99,35){\circle*{0.61}}
\put(6,35.5){\line(-6,4){6}}
\put(27,35.5){\line(-4,4){4}}
\put(44,35.5){\line(0,1){4}}
\put(58,35.5){\line(0,1){4}}
\put(76,35.5){\line(0,1){4}}
\put(91,35.5){\line(0,1){4}}
\put(0,40){\circle*{0.61}}
\put(9,40){\circle*{0.61}}
\put(23,40){\circle*{0.61}}
\put(35,40){\circle*{0.61}}
\put(44,40){\circle*{0.61}}
\put(50,40){\circle*{0.61}}
\put(58,40){\circle*{0.61}}
\put(69,40){\circle*{0.61}}
\put(76,40){\circle*{0.61}}
\put(84,40){\circle*{0.61}}
\put(91,40){\circle*{0.61}}
\put(99,40){\circle*{0.61}}
\put(23,40.5){\line(-10,4){10}}
\put(44,40.5){\line(-13,4){13}}
\put(58,40.5){\line(-9,4){9}}
\put(76,40.5){\line(0,1){4}}
\put(91,40.5){\line(0,1){4}}
\put(0,45){\circle*{0.61}}
\put(13,45){\circle*{0.61}}
\put(26,45){\circle*{0.61}}
\put(31,45){\circle*{0.61}}
\put(41,45){\circle*{0.61}}
\put(49,45){\circle*{0.61}}
\put(61,45){\circle*{0.61}}
\put(76,45){\circle*{0.61}}
\put(84,45){\circle*{0.61}}
\put(91,45){\circle*{0.61}}
\put(99,45){\circle*{0.61}}
\put(13,45.5){\line(-13,4){13}}
\put(31,45.5){\line(-7,4){7}}
\put(49,45.5){\line(0,1){4}}
\put(76,45.5){\line(0,1){4}}
\put(91,45.5){\line(0,1){4}}
\put(0,50){\circle*{0.61}}
\put(14,50){\circle*{0.61}}
\put(24,50){\circle*{0.61}}
\put(37,50){\circle*{0.61}}
\put(49,50){\circle*{0.61}}
\put(61,50){\circle*{0.61}}
\put(76,50){\circle*{0.61}}
\put(84,50){\circle*{0.61}}
\put(91,50){\circle*{0.61}}
\put(99,50){\circle*{0.61}}
\put(24,50.5){\line(-16,4){16}}
\put(49,50.5){\line(-19,4){19}}
\put(76,50.5){\line(-29,4){29}}
\put(91,50.5){\line(-24,4){24}}
\put(81,55){\line(6,-1){19}}
\put(96,55){\line(6,-1){4}}
\put(0,55){\circle*{0.61}}
\put(8,55){\circle*{0.61}}
\put(19,55){\circle*{0.61}}
\put(30,55){\circle*{0.61}}
\put(36,55){\circle*{0.61}}
\put(47,55){\circle*{0.61}}
\put(59,55){\circle*{0.61}}
\put(67,55){\circle*{0.61}}
\put(73,55){\circle*{0.61}}
\put(81,55){\circle*{0.61}}
\put(88,55){\circle*{0.61}}
\put(96,55){\circle*{0.61}}
\put(-13,21.5){time $t$}
\put(-3.6,-1){0}
\put(-3.6,4){1}
\put(-3.6,9){2}
\put(-3.6,14){3}
\put(-3.6,19){4}
\put(-3.6,24){5}
\put(-3.6,29){6}
\put(-3.6,34){7}
\put(-3.6,39){8}
\put(-3.6,44){9}
\put(-4.6,49){10}
\put(-4.6,54){11}
\put(0,-2.5){\line(1,0){100}}
\multiput(0,-2.5)(10,0){11}{\line(0,-1){1.3}}
\put(41,-10){position $x$}
\put(-0.8,-7){0}
\put(19.2,-7){2}
\put(39.2,-7){4}
\put(59.2,-7){6}
\put(79.2,-7){8}
\put(98.2,-7){10}
\put(103,10){$\bullet$}
\end{picture}

\vspace{0.7in}

{\bf Figure 1.} 
A realization of the process, showing space-time trajectories of alternate customers.  
The customer with initial rank 12 is marked by $\bullet$.

\medskip
Figure 1 shows a realization of the process over the time interval
$0 \le t \le 11$ and the spatial range 
$0 \le x \le 10$.  
Time increases upwards, and the head of the queue is on the left. 
The $\bullet$ indicate customer positions at each time, and
the lines indicate the space-time trajectories of individual customers 
(for visual clarity, only alternate customers' trajectories are shown).

Although Figure 1 seems the intuitively natural way to draw a realization, a different graphic is more  suggestive for mathematical analysis.  
A configuration 
$\bx = (0 = x_0 < x_1 < x_2 < x_3 \ldots)$ 
can be represented by its centered counting function 
\begin{equation}
F(x) := \max\{k: \  x_k \le x\} \ - x , \quad 0 \le x < \infty .
\label{def-ccc}
\end{equation}
as illustrated in Figure 2.

\begin{picture}(100,20)(0,50)
\put(0,55){\circle*{0.91}}
\put(14,55){\circle*{0.91}}
\put(24,55){\circle*{0.91}}
\put(30,55){\circle*{0.91}}
\put(36,55){\circle*{0.91}}
\put(47,55){\circle*{0.91}}
\put(59,55){\circle*{0.91}}
\put(67,55){\circle*{0.91}}
\put(73,55){\circle*{0.91}}
\put(81,55){\circle*{0.91}}
\put(88,55){\circle*{0.91}}
\put(96,55){\circle*{0.91}}
\put(0,55){\line(2,-1){14}}
\put(14,48){\line(0,1){5}}
\put(14,53){\line(2,-1){10}}
\put(24,48){\line(0,1){5}}
\put(24,53){\line(2,-1){6}}
\put(30,50){\line(0,1){5}}
\put(30,55){\line(2,-1){6}}
\put(36,52){\line(0,1){5}}
\put(36,57){\line(2,-1){11}}
\put(47,51.5){\line(0,1){5}}
\put(47,56.5){\line(2,-1){12}}
\put(59,50.5){\line(0,1){5}}
\put(59,55.5){\line(2,-1){8}}
\put(67,51.5){\line(0,1){5}}
\put(67,56.5){\line(2,-1){6}}
\put(73,53.5){\line(0,1){5}}
\put(73,58.5){\line(2,-1){8}}
\put(81,54.5){\line(0,1){5}}
\put(81,59.5){\line(2,-1){7}}
\put(88,56){\line(0,1){5}}
\put(88,61){\line(2,-1){8}}
\put(96,57){\line(0,1){5}}
\put(96,62){\line(2,-1){4}}
\put(0,47.5){\line(1,0){100}}
\multiput(0,47.5)(10,0){11}{\line(0,-1){1.3}}
\put(41,40){position $x$}
\put(19.2,43){2}
\put(39.2,43){4}
\put(59.2,43){6}
\put(79.2,43){8}
\put(98.2,43){10}
\put(0,43){\line(0,1){27}}
\put(0,50){\line(-1,0){1.3}}
\put(0,55){\line(-1,0){1.3}}
\put(0,60){\line(-1,0){1.3}}
\put(-3.9,54){0}
\put(-3.9,59){1}
\put(-4.9,49){-1}
\end{picture}

\vspace{0.5in}

{\bf Figure 2.} 
The centered counting function associated with a configuration 
$\bullet$.

\medskip
At each time $t$, let us consider the centered counting function 
$F_t(x)$ and plot the graph of the upward-translated function
\begin{equation}
x \to G(t,x) := t + F_t(x) .
\label{def-G}
\end{equation}
In other words, we draw the function starting at the point $(0,t)$ instead of the origin.  
Taking the same realization as in Figure 1, and
superimposing all these graphs, 
gives Figure 3.

\newpage

\begin{picture}(100,65)
\put(0,0){\line(2,-1){12}}
\put(12,-6){\line(0,1){5}}
\put(12,-1){\line(2,-1){10}}
\put(22,-6){\line(0,1){5}}
\put(22,-1){\line(2,-1){7}}
\put(29,-4.5){\line(0,1){5}}
\put(29,0.5){\line(2,-1){7}}
\put(36,-3){\line(0,1){5}}
\put(36,2){\line(2,-1){14}}
\put(50,-5){\line(0,1){5}}
\put(50,0){\line(2,-1){13}}
\put(63,-6.5){\line(0,1){5}}
\put(63,-1.5){\line(2,-1){7}}
\put(70,-5){\line(0,1){5}}
\put(70,0){\line(2,-1){9}}
\put(79,-4.5){\line(0,1){5}}
\put(79,0.5){\line(2,-1){7}}
\put(86,-3){\line(0,1){5}}
\put(86,2){\line(2,-1){8}}
\put(94,-2){\line(0,1){5}}
\put(94,3){\line(2,-1){6}}
\put(0,5){\line(2,-1){7}}
\put(7,1.5){\line(0,1){5}}
\put(7,6.5){\line(2,-1){14}}
\put(21,-0.5){\line(0,1){5}}
\put(21,4.5){\line(2,-1){15}}
\put(0,10){\line(2,-1){10}}
\put(10,5){\line(0,1){5}}
\put(10,10){\line(2,-1){10}}
\put(20,5){\line(0,1){5}}
\put(20,10){\line(2,-1){15}}
\put(35,2.5){\line(0,1){5}}
\put(35,7.5){\line(2,-1){9}}
\put(44,3){\line(0,1){5}}
\put(44,8){\line(2,-1){16}}
\put(60,0){\line(0,1){5}}
\put(60,5){\line(2,-1){12}}
\put(72,-1){\line(0,1){5}}
\put(72,4){\line(2,-1){14}}
\put(0,15){\line(2,-1){7}}
\put(7,11.5){\line(0,1){5}}
\put(7,16.5){\line(2,-1){10}}
\put(17,11.5){\line(0,1){5}}
\put(17,16.5){\line(2,-1){8}}
\put(25,12.5){\line(0,1){5}}
\put(25,17.5){\line(2,-1){9}}
\put(34,13){\line(0,1){5}}
\put(34,18){\line(2,-1){12}}
\put(46,12){\line(0,1){5}}
\put(46,17){\line(2,-1){9}}
\put(55,12.5){\line(0,1){5}}
\put(55,17.5){\line(2,-1){6}}
\put(61,14.5){\line(0,1){5}}
\put(61,19.5){\line(2,-1){8}}
\put(69,15.5){\line(0,1){5}}
\put(69,20.5){\line(2,-1){9}}
\put(78,16){\line(0,1){5}}
\put(78,21){\line(2,-1){13}}
\put(78,21){\circle*{1.5}}
\put(91,14.5){\line(0,1){5}}
\put(91,19.5){\line(2,-1){9}}
\put(100,15){\line(0,1){5}}
\put(0,20){\line(2,-1){13}}
\put(13,13.5){\line(0,1){5}}
\put(13,18.5){\line(2,-1){12}}
\put(0,25){\line(2,-1){5}}
\put(5,22.5){\line(0,1){5}}
\put(5,27.5){\line(2,-1){8}}
\put(13,23.5){\line(0,1){5}}
\put(13,28.5){\line(2,-1){7}}
\put(20,25){\line(0,1){5}}
\put(20,30){\line(2,-1){13}}
\put(33,23.5){\line(0,1){5}}
\put(33,28.5){\line(2,-1){14}}
\put(47,21.5){\line(0,1){5}}
\put(47,26.5){\line(2,-1){9}}
\put(56,22){\line(0,1){5}}
\put(56,27){\line(2,-1){5}}
\put(61,24.5){\line(0,1){5}}
\put(61,29.5){\line(2,-1){13}}
\put(61,29,5){\circle*{1.5}}
\put(74,23){\line(0,1){5}}
\put(74,28){\line(2,-1){6}}
\put(80,25){\line(0,1){5}}
\put(80,30){\line(2,-1){7}}
\put(87,26.5){\line(0,1){5}}
\put(87,31.5){\line(2,-1){10}}
\put(97,26.5){\line(0,1){5}}
\put(97,31.5){\line(2,-1){3}}
\put(0,30){\line(2,-1){10}}
\put(10,25){\line(0,1){5}}
\put(10,30){\line(2,-1){10}}
\put(20,25){\line(0,1){5}}
\put(0,35){\line(2,-1){6}}
\put(6,32){\line(0,1){5}}
\put(6,37){\line(2,-1){11}}
\put(17,31.5){\line(0,1){5}}
\put(17,36.5){\line(2,-1){10}}
\put(27,31.5){\line(0,1){5}}
\put(27,36.5){\line(2,-1){8}}
\put(35,32.5){\line(0,1){5}}
\put(35,37.5){\line(2,-1){9}}
\put(44,33){\line(0,1){5}}
\put(44,38){\line(2,-1){6}}
\put(44,38){\circle*{1.5}}
\put(50,35){\line(0,1){5}}
\put(50,40){\line(2,-1){8}}
\put(58,36){\line(0,1){5}}
\put(58,41){\line(2,-1){11}}
\put(69,35.5){\line(0,1){5}}
\put(69,40.5){\line(2,-1){7}}
\put(76,37){\line(0,1){5}}
\put(76,42){\line(2,-1){8}}
\put(84,38){\line(0,1){5}}
\put(84,43){\line(2,-1){7}}
\put(91,39.5){\line(0,1){5}}
\put(91,44.5){\line(2,-1){8}}
\put(99,40.5){\line(0,1){5}}
\put(99,45.5){\line(2,-1){1}}
\put(0,40){\line(2,-1){9}}
\put(9,35.5){\line(0,1){5}}
\put(9,40.5){\line(2,-1){14}}
\put(23,33.5){\line(0,1){5}}
\put(23,38.5){\line(2,-1){12}}
\put(35,32.5){\line(0,1){5}}
\put(0,45){\line(2,-1){13}}
\put(13,38.5){\line(0,1){5}}
\put(13,43.5){\line(2,-1){13}}
\put(26,37){\line(0,1){5}}
\put(26,42){\line(2,-1){5}}
\put(31,39.5){\line(0,1){5}}
\put(31,44.5){\line(2,-1){10}}
\put(31,44.5){\circle*{1.5}}
\put(41,39.5){\line(0,1){5}}
\put(41,44.5){\line(2,-1){8}}
\put(49,40.5){\line(0,1){5}}
\put(49,45.5){\line(2,-1){12}}
\put(61,39.5){\line(0,1){5}}
\put(61,44.5){\line(2,-1){15}}
\put(76,37){\line(0,1){5}}
\put(0,50){\line(2,-1){14}}
\put(14,43){\line(0,1){5}}
\put(14,48){\line(2,-1){10}}
\put(24,43){\line(0,1){5}}
\put(24,48){\line(2,-1){13}}
\put(24,48){\circle*{1.5}}
\put(37,41.5){\line(0,1){5}}
\put(37,46.5){\line(2,-1){12}}
\put(49,40.5){\line(0,1){5}}
\put(0,55){\line(2,-1){8}}
\put(8,51){\line(0,1){5}}
\put(8,56){\line(2,-1){11}}
\put(8,56){\circle*{1.5}}
\put(19,50.5){\line(0,1){5}}
\put(19,55.5){\line(2,-1){11}}
\put(30,50){\line(0,1){5}}
\put(30,55){\line(2,-1){6}}
\put(36,52){\line(0,1){5}}
\put(36,57){\line(2,-1){11}}
\put(47,51.5){\line(0,1){5}}
\put(47,56.5){\line(2,-1){12}}
\put(59,50.5){\line(0,1){5}}
\put(59,55.5){\line(2,-1){8}}
\put(67,51.5){\line(0,1){5}}
\put(67,56.5){\line(2,-1){6}}
\put(73,53.5){\line(0,1){5}}
\put(73,58.5){\line(2,-1){8}}
\put(81,54.5){\line(0,1){5}}
\put(81,59.5){\line(2,-1){7}}
\put(88,56){\line(0,1){5}}
\put(88,61){\line(2,-1){8}}
\put(96,57){\line(0,1){5}}
\put(96,62){\line(2,-1){4}}
\put(-13,21.5){time $t$}
\put(-3.6,-1){0}
\put(-3.6,4){1}
\put(-3.6,9){2}
\put(-3.6,14){3}
\put(-3.6,19){4}
\put(-3.6,24){5}
\put(-3.6,29){6}
\put(-3.6,34){7}
\put(-3.6,39){8}
\put(-3.6,44){9}
\put(-4.6,49){10}
\put(-4.6,54){11}
\put(0,-12.5){\line(1,0){100}}
\multiput(0,-12.5)(10,0){11}{\line(0,-1){1.3}}
\put(41,-20){position $x$}
\put(-0.8,-17){0}
\put(19.2,-17){2}
\put(39.2,-17){4}
\put(59.2,-17){6}
\put(79.2,-17){8}
\put(98.2,-17){10}

\end{picture}

\vspace{1.1in}

{\bf Figure 3.}
A realization of the process $G(t,x)$.  The $\bullet$ indicate the points associated with the moves of the customer with initial rank 12; compare with the space-time trajectory in Figure 1.

\medskip
Clearly the graphs in Figure 3 are ``coalescing" in some sense.
What is precisely true is the following. The assertion
\begin{quote}
 the rank $k$ person (for some $k$)  at time $t$ is at position $x^*$ and
remains at that position $x^*$ at time $t+1$
\end{quote}
is equivalent to the assertion
\begin{quote}
the graphs starting at $(0,t)$ and at $(0,t+1)$ both contain the same step 
of the form $(x^*,s)$ to $(x^*,s+1)$ (for some $s$).
\end{quote}
The former assertion implies that customers at ranks greater than $k$ also do not move, and therefore the graphs in the second assertion 
coincide for all $x \ge x^*$.
So, even though the graphs may not coalesce upon first meeting 
(as is clear from the left region of Figure 3), they do coalesce at the first point where they make the same jump.

Many readers will immediately perceive that Figure 3 shows 
``something like" coalescing random walks in one dimension.  
Our goal is to formalize this as an asymptotic theorem.
But note that, for a graphic of a model of coalescing random walks which 
is drawn to  resemble Figure 3, the horizontal axis would be ``time" and the vertical axis would be ``space".  So in making the analogy between our process and 
coalescing random walks, we are interchanging the roles of space and time. 
In particular, in Figure 3 our process is (by definition) 
Markov in the vertical direction, whereas models of coalescing random walks pictured in Figure 4 below are  (by definition) Markov in the horizontal direction.  
These means that  in conditioning arguments we will need to be careful, in fact somewhat obsessive,  to be clear about what exactly is being conditioned upon.

\subsection{Coalescing Brownian motion on the line}
\label{sec:CBM}
Given $\tau > 0$ and a random locally finite subset $\eta_\tau \subset \Reals$, it is straightforward to define a process 
$(B_y(t), \tau \le t < \infty, y \in \eta_\tau)$ 
which formalizes the idea

\smallskip
\noindent
(i) at time $\tau$, start independent Brownian motions $(B^0_y(t), t \ge \tau)$ 
on $\Reals$  from  each $y \in \eta_\tau$;

\noindent
(ii)  when two of these diffusing ``particles" meet, they coalesce into one particle which continues diffusing as Brownian motion.

\smallskip
\noindent
So $B_y(t)$ is the time-$t$ position of the particle containing the particle started at $y$.
Call this process the  CBM$(\eta_\tau, \tau)$ process.
A line of research starting with Arratia \cite{arratia1979}  (see also \cite{toth-werner}) 
shows this construction can be extended to the {\em standard CBM (coalescing Brownian motion) process}
$(B(y,t), - \infty < y < \infty, \  t \ge 0)$  
which is determined by the properties 

\smallskip
\noindent (iiii) 
given the point process $\eta_\tau = \{B(y,\tau), \ y \in \Reals\}$ of time-$\tau$ particle positions, the process subsequently
evolves as CBM$(\eta_\tau, \tau)$ 

\noindent (iv)
$\lim_{\tau \downarrow 0} \Pr( \eta_\tau \cap I = \emptyset) = 0  \quad \mbox{ for each open } I \subset \Reals $.

\smallskip \noindent
The latter requirement formalizes the idea that the process starts at time $0$ with a particle at {\em every} point 
$y \in \Reals$.  
In particular,  the point process 
$\eta_t = \{B(y, t), \ y \in \Reals\}$ of time-$t$ particle positions
inherits the scaling and translation-invariance properties of Brownian motion:

\smallskip
\noindent
(v) For each $t>0$, the point process $\eta_t(\cdot)$ is translation-invariant and so has some mean density $\rho_t$.

\noindent
(vi) The distribution of the family  $(\eta_t(\cdot), 0<t<\infty)$ is invariant under the scaling map 
$(t,y) \to (ct,c^{1/2}y)$.

\noindent
(vii) The distribution of $\eta_t(\cdot)$ is the push-forward of the distribution of $\eta_1(\cdot)$ under
the map $y \to t^{1/2}y$ and so 
\begin{equation}
\rho_t = t^{-1/2} \rho_1 .
\label{def-rho}
\end{equation}

\medskip
What will be important in this paper is an ``embedding property". 
 If we take the CBM$(\tilde{\eta}^0,0)$ process, starting at time $0$ 
with some translation-invariant 
locally finite point process $\tilde{\eta}^0$ with mean density $0 < \tilde{\rho}_0 < \infty$, and 
embed this process into standard CBM in the natural way, 
it is straightforward to see that the time-asymptotics of CBM$(\tilde{\eta}^0,0)$ are the same as the 
time-asymptotics of standard CBM.

\subsection{What does the CBM limit tell us?}
\label{sec:what}
Figure 3 suggests an asymptotic result of the form
\begin{quote}
the spatial queue process, in the Figure 3 representation as a process $G(t,x)$, converges in distribution 
(in the usual random walk to Brownian motion scaling limit) to the CBM process.
\end{quote}
We formalize this in section \ref{sec:CBMlimit} as Theorem \ref{T:main}, by working on
 a  ``coalescing continuous functions" state space.
Our formalization requires a double limit:  we study the process as it affects  customers beyond position $x$, in a  time $t \to \infty$ limit,
and show the $x \to \infty$ rescaled limit is a CBM process.  
And as in our earlier statement of Theorem \ref{T:main_rate},
we can only prove this for  time-averaged limits.

Now the Figure 3 representation as a process $G(t,x)$ is a somewhat ``coded" version of the more natural
Figure 1 graphic of customers' trajectories, so to interpret the CBM scaling limit 
we need to ``decode" Figure 3.
The wave lengths, measured by position of first unmoved customer, 
appear in Figure 3 as the positions $x^\prime$ where two functions $G(t^\prime,\cdot), G(t^{\prime \prime},\cdot)$ 
coalesce.
In other words the time-asymptotic rate of waves of length $>x$ equals the time-asymptotic rate at which distinct functions
$G(t,\cdot)$ cross over $x$.
In the space-asymptotic limit (of rank $j$ or position $x$),
as a consequence  of results used to prove the 
order of magnitude bounds (Proposition \ref{P-main}) we have
that $x/j \to 1$.
Now the  scaling limit identifies ``distinct functions
$G(t,\cdot)$ crossing over $x$" with CBM particles at time $x$,
 and so Theorem \ref{T:main} will imply a result for the
first-order asymptotics of wave length, measured by either rank or position, which sharpens Proposition \ref{P-main}. 
Recall $\sigma$ is the s.d. of $\mu$.
\begin{Theorem}
\label{T:main_rate}
For any initial distribution, 
\[ \limsup_{\tau \to \infty} \tau^{-1} \sum_{t=1}^\tau \Pr(W(t) > i) \sim \rho_1 \sigma^{-1} i^{-1/2} \mbox{ as }i \to \infty \]
\[ \liminf_{\tau \to \infty} \tau^{-1} \sum_{t=1}^\tau \Pr(W(t) > i) \sim \rho_1 \sigma^{-1} i^{-1/2} \mbox{ as }i \to \infty \]
where $\rho_1$ is the numerical constant (\ref{def-rho}) associated with CBM. 
\end{Theorem}
This will be a consequence of Corollary \ref{C:rateU}.
As discussed earlier, if we knew there was convergence to stationarity we could restate this in the simpler way
\begin{equation}
\mbox{at stationarity, $\Pr(W(t) > i) \sim \rho_1 \sigma^{-1} i^{-1/2}$ as $i \to \infty$}.
\label{single-limit}
\end{equation}

\vspace*{-0.1in}

\includegraphics[width=4.5in]{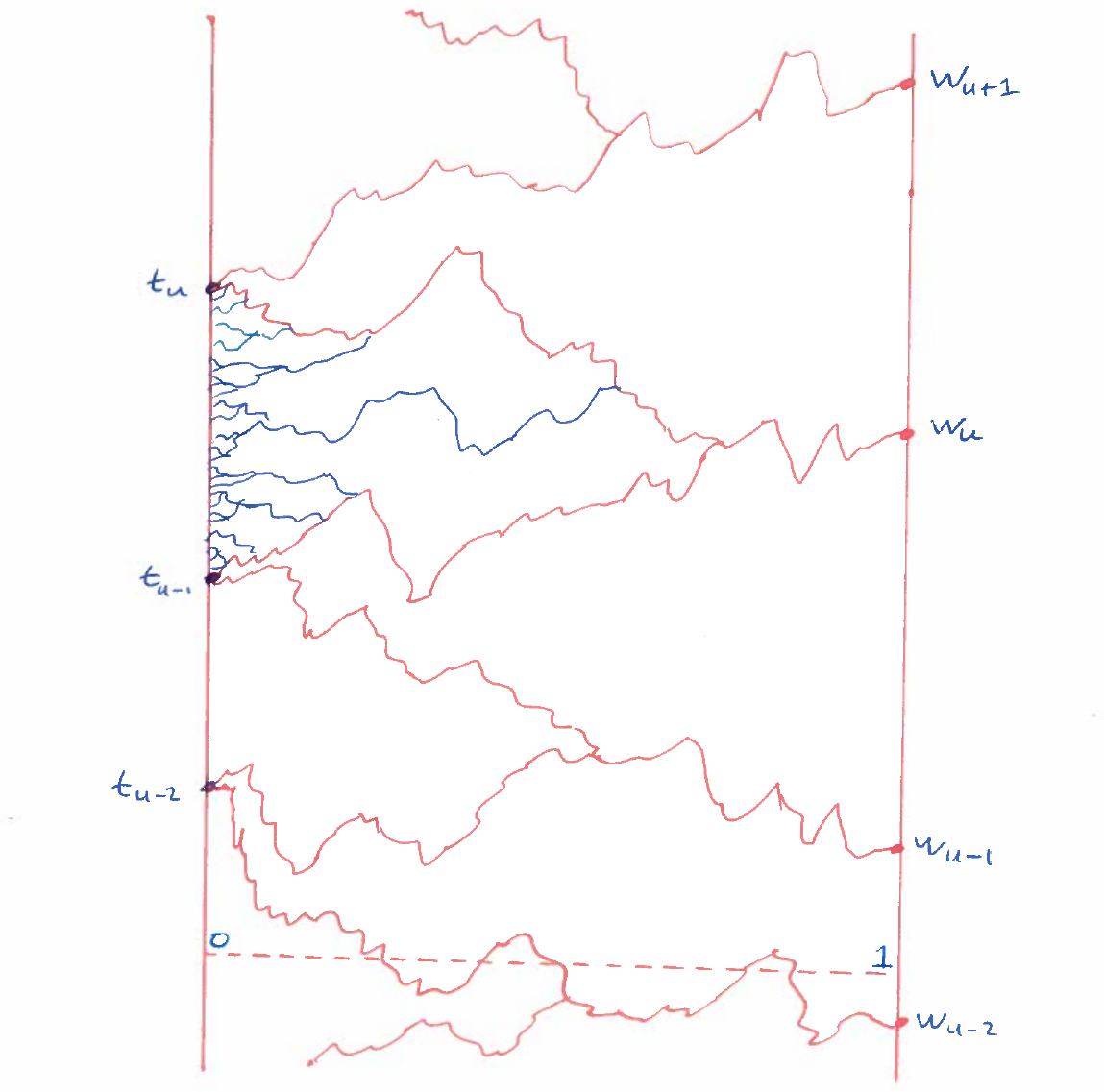}

\vspace*{-0.9in}

{\bf Figure 4.}
Sketch of CBM over the time interval $[0,1]$.

\bigskip
The implications of the CBM scaling limit for the space-time trajectory of a typical customer are rather
more subtle.
As suggested by Figure 8, the trajectory is deterministic to first order, so the CBM limit must be telling us something about 
the second-order behavior.

Look at Figure 4, which shows 
CBM over ``time" $0 \le x \le 1$ and ``space" $- \infty < t < \infty$.
There are positions $(w_u, \ u \in \Ints)$ of the particles at time $x = 1$;  and
\begin{quote}
(*) 
each such particle is the coalescence of the particles which started at $x=0$ in some interval $(t_{u-1},t_u)$.
\end{quote}
So we have two stationary point processes on $\Reals$, the 
processes $(w_u, \ u \in \Ints)$ and $(t_u, \ u \in \Ints)$, which must have the same rate $\rho_1$.

In the queue model, as Figure 3 illustrates, the motion of one particular individual can be seen 
within the family of graphs as line segments within one diagonal line, To
 elaborate, write
$(x_u,s_u)$ for the coordinates of the top end points (the $\bullet$'s in Figure 3) of these segments, after the $u$Õth move of the individual.
The $u$Õth move happens at some time $t_u$, when the individual moves to position $x_u$, and so 
$s_u = G(t_u,x_u)$.
Until the next move, that is for $t_u \le t < t_{u+1}$, 
we have (from the definition of $G$) that $G(t,x_u) = s_u$.  
At the time $t_{u+1}$ of the next move, the individual moves to some position $x_{u+1} < x_u$. 
From the definition of $G$ we have 
\[ s_{u+1} - s_u = G(t_{u+1},x_{u+1} ) - G(t_u,x_u) = - x_{u+1} + x_u .\]
So the points $(x_u,s_u)$ lie on a diagonal line
\[ x_u + s_u = \mbox{ constant} . \]
In fact the constant is the initial rank of the individual, 
or equivalently the time at which the individual will reach the head of the queue. 
What is important to note for our analysis is that the times $t_u$ of the moves are not 
indicated directly by the points $(x_u,s_u)$, but instead are specified by 
\begin{equation}
G(t,x_u) = s_u, \ t_u \le t < t_{u+1} .
\end{equation}
In words, 
\begin{quote}
(**)  the $t_u$ are the end-times of the intervals of $t$ for which the functions $G(t,\cdot)$ have coalesced before the position of the individual.
\end{quote}

So consider a customer at position approximately  $n$, and hence with rank approximately $n$, for large $n$.
From the order-of-magnitude results for wave lengths, the time between moves, and the distances moved, 
for this customer are both order $n^{1/2}$.
But in the Brownian rescaling which takes  Figure 3  to Figure 4, the horizontal axis is rescaled by $1/n$ and the vertical axis by 
$1/n^{1/2}$.  That is, 1 unit time in the CBM model corresponds to $n$ space units in the queue model, and 1 space unit in the CBM model corresponds to $n^{1/2}$ time units in the queue model.  
So when our position-$n$ customer moves, in the CBM approximation the horizontal distance moved is only
$O(n^{-1/2})$  whereas the vertical distance is order $1$. 
So in the scaling  limit the Figure 3 points corresponding to that customer become 
 the points $w_u$  on the vertical line in Figure 4.
One consequence of this limit is provided by  (*) and (**): 
the sequence of times of moves of that customer becomes, in the scaling limit, the sequence $(t_u)$ 
 in the CBM process indicating the spatial intervals of the time-$0$ particles that have coalesced into the time-$1$ particles.
Precisely, from Theorem \ref{T:main} we can deduce (section \ref{sec:s-t} gives the outline) the following corollary, 
in which the first component formalizes the assertion above.

Let $U_\tau$ have uniform distribution on $\{1,2,\ldots,\tau\}$.
\begin{Corollary}
\label{C:2C}
In the queue model, consider the individual at time $U_{\tau_n}$ nearest to position $n$.
Write $((t^{(n)}_u, w^{(n)}_u), \ u \ge 1)$ for the times $t$ and end-positions $w$ of moves of that individual.  
If $\tau_n \to \infty$ sufficiently fast then
\[
\left( \frac{t^{(n)}_u - U_{\tau_n}}{n^{1/2}} \ , \ \frac{n - w^{(n)}_u }{n^{1/2}} \right)_{u \ge 1}
\cd
 \big(\sigma t_u, \sigma w_u \big)_{- \infty < u < \infty}
\]
in the sense of convergence of point processes on $\Reals^2$.
\end{Corollary}
If we knew there was convergence to stationarity in the queue model then we could use $\tau_n$ in place of $U_{\tau_n}$.

How does the second component in Corollary \ref{C:2C} arise?  
Return to the setting of (**), and consider the graphs 
$G(t_{u+1}, \cdot)$ and $G(t_u,\cdot)$ 
coding  configurations at successive moves of the individual under consideration.
From the definition of $G$, the increment 
$G(t_{u+1}, x_u) - G(t_u, x_u)$ equals (up to $\pm O(1)$) the number ($m$, say) of customers who cross position $x_u$ in the time-$(t_{u+1})$ wave.
So the distance moved by the individual under consideration. is essentially the sum of these $m$ inter-customer spacings. 
These spacings are  i.i.d.($\mu$) subject to a certain conditioning, but the conditioning event does not have small probability, so the distance moved by that customer is $(1 \pm o(1))m$.  
The size of an increment  $G(t_{u+1}, x) - G(t_u, x)$ as a function of $x$ is
first-order constant over the spatial interval $x \in (x_u \pm O(x_u^{1/2}))$ where the individual will be for the next $O(1)$ moves.
So for large $n$ and $t$, the successive distances moves by a rank $n$ customer around time $t$ are to first order the 
distances between the successive distinct values of $t \to G(t,n)$.
In the scaling limit these become the successive increments of $(w_u)$ (Figure 4) in the CBM process, 
explaining the second component in Corollary \ref{C:2C}.

We will add some details to this outline  in section \ref{sec:s-t}.

\subsection{Outline of proofs}
\label{sec:outline}
We found it surprisingly hard to devise proofs, in part because of the subtleties, mentioned earlier, arising in formalizing conditioning arguments.
Simpler proofs may well be discovered in future.

The main goal of this paper is to prove Proposition \ref{P-main}, that waves of lengths $>k$ occur at rate of order $k^{-1/2}$, not faster or slower.  
As one methodology, we can analyze the distribution at a large time by looking backwards, because 
the current configuration must consist of successive ``blocks" of customers, each block having moved to their 
current position at the same past time.  In particular,
because the inter-customer distances within a block are roughly i.i.d.($\mu$) (only {\em roughly} because of conditioning effects),
and because block lengths tend to grow, we expect that 
(in our conjectured stationary limit) the space-asymptotic inter-customer distances 
should become approximately i.i.d.($\mu$).  
One formalization of this idea is given in 
Lemma \ref{L:almost}, and this leads (section \ref{sec:lb}) to an order $k^{-1/2}$ lower bound on the rate of
waves of lengths $>k$; but for technical reasons we first need to know the corresponding upper bound.
Our proof of the upper bound (outlined at the start of section \ref{sec:UB}) 
is intricate, but the key starting point is the observation (section \ref{sec:close_waves}) that, 
after a long wave, the inter-customer spacings of the moved customers tend to be shorter than i.i.d.($\mu$); 
this and classical random walk estimates lead to an order $k^{1/2}$ lower bound on the mean time between successive 
waves of lengths $>k$, which is an  order $k^{-1/2}$ upper bound on the rate of such waves.
Although intricate in detail, all these proofs essentially rely only on explicit bounds, in contrast to those mentioned below.

In section \ref{sec:CBMlimit} we develop the CBM scaling limit results, and their implications, discussed 
in section \ref{sec:what}.
Given the order of magnitude result (Proposition \ref{P-main}), this further development is essentially ``soft", involving weak convergence 
and compactness on an unfamiliar space (of coalescing continuous functions), combined with basic properties of CBM 
and classical weak convergence of random walks to Brownian motion.
We give a careful outline of such an analysis without attending to all the technical details.

\subsection{Remarks on the model and related work} 
Conceptually the model 
seems fairly realistic and to provide a plausible explanation for 
what the reader might have observed in long queues such as airport security lines. 
See section \ref{sec:OAK} for some data. 
In reality there are other reasons why a customer might not move when the `wave" reaches them -- not paying
attention because of talking or texting, for instance. 
Two  assumptions of our model --  that service times are constant, and that the successive customers 
in a wave move simultaneously -- are 
unrealistic, but these are irrelevant to the particular feature (length of waves of motion) that we are studying.
Intuitively the model behavior should also be robust to different customers having different parameters 
$\mu, c_-,c^+$ -- see section \ref{sec:robust}.

We do not know any previous work on closely related models.
In one way, our model could be regarded as a 
continuum, syncronous-moves  analog of the (discrete-space, asynchronous moves) TASEP which has been intensively studied \cite{borodin}
owing to its connections with other statistical physics models.  The small literature on such analogs of TASEP (see \cite{blank_2010}  and citations therein) has focused on ergodic properties for processes on the doubly-infinite line.
In another way,  our phenomenon could be regarded as analogous to stop-and-go motion in traffic jams, 
for which realistic models would need to consider both positions and velocities of vehicles, although discrete-space
discrete-time synchronous models have also been studied  \cite{gray}.
Academic literature such as  \cite{doi:10.3141/2177-08} studying real world security lines does not 
address the specific ``wave" feature studied in this paper.

Finally, our concern in this paper is different from the usual concerns  of {\em Stochastic Systems} 
(evaluation or optimization of some measure of system performance).  
Instead the paper seeks to articulate a plausible probabilistic explanation of an observed everyday phenomenon,  
in the spirit of the author's "probability in the real world" project \cite{me-RW}.

\newpage
\section{Upper bounding the rate of long waves}
\label{sec:UB}
In this section we show that the rate of waves of length $>j$ is $O(j^{-1/2})$.
Informally, we do this by considering the progress of an individual from rank $j + o(j)$ to rank $j - o(j)$.
The individual only moves when there is a wave of length $> j \pm o(j)$.  
The key idea is 
\begin{quote}
If the individual moves a distance $o(j^{1/2})$ then 
the next wave that moves the individual  is unlikely to have length $> 2j$.
\end{quote}
The argument is intricate and in fact uses that idea somewhat indirectly (Proposition \ref{P:Xst}).
The remaining argument is outlined in section \ref{sec:ub}.

Figure 3 seems the most intuitive way to picture waves, in terms of their spatial extent, but for proofs it is
more convenient to work with ranks.
To track the progress of an individual we use the notation
\[
X^*_s(t) = X_{s-t}(t), \ 0 \le t \le s . \]
That is, $X^*_s(\cdot)$ tracks the position of ``individual $s$", the individual who starts at rank s and will reach the service position at time $s$.
 Write 
 \[ \bX_{[s]}(t) = (X_i(t), \ 0 \le i \le s-t) \]
 for the positions of ``individual $s$"  and the customers ahead of that individual in the queue.
Write
\[ \FF(t) = \sigma(\bX(u), 0 \le u \le t) \]
for the natural filtration of our process, which is constructed in terms of i.i.d($\mu$) random variables 
$\{\xi_i(t), i \ge1, t \ge 1\}$
as at (\ref{Xit1} -  \ref{Xit3}).

\subsection{Close waves coalesce quickly}
\label{sec:close_waves}
The event that individual $s$ is moved in the time-$t$ wave is the event
\[ 
\{X^*_s(t) < X^*_s(t-1) \} = \{W(t) > s-t\} .
\]
Such a wave may be a ``long" wave, in that it extends at least $j$ customers past individual $s$.
Our first goal is to show, roughly speaking, that after one long wave, and before individual $s$ 
moves a distance $o(j^{1/2})$:

\medskip \noindent
either (i)  there is unlikely to be another long wave; \\
or (ii) there are likely to be many waves that move individual $s$.

\medskip \noindent
This is formalized by 
Propositions \ref{P:Xst} and \ref{P:NCA}, though 
complicated by the fact we will need to consider ``good" long waves.

How far one wave reaches depends on the growth rate 
of the partial sums $\sum_{i=1}^k \xi_i$ associated with this new wave, relative to the partial sums 
associated with the previous long wave. 
Because we are interested in the difference, 
what is ultimately relevant is the distribution of the maximum of the symmetrized (hence mean zero) random walk
\begin{equation}
S^{sym}_k = \sum_{i=1}^k (\xi^\prime_i - \xi^{\prime\prime}_i)
\label{Ssym}
\end{equation}
where $(\xi^\prime_i)$ and $(\xi^{\prime \prime}_i)$ denote
 independent i.i.d.($\mu)$  sequences. 
Define 
\begin{equation}
q(j,y) =  \Pr(\max_{1 \le k \le j} S^{sym}_k \le y ) .
\label{def-qjy}  
\end{equation}
For later use, note that  Donsker's theorem and the fact that the distribution of the maximum
of  Brownian motion over $0 \le t \le 1$ has bounded  density imply
\begin{Lemma}
\label{LRW}
There exists a constant $C$ and $\delta_j \downarrow 0$, depending only on the distribution $\mu$, such that
\[ q(j,y) \le Cy j^{-1/2} + \delta_j, \ \ y \ge 1, j \ge 1 .\]
\end{Lemma}

Now consider
\begin{equation}
q_{j,y}(x_1,\ldots,x_j) =  \Pr \left( \max_{1 \le k \le j}   \sum_{i=1}^k (\xi^\prime_i - x_i)   \le y \right)   .
\label{def-qjyx}
\end{equation}
Say a sequence $(x_1,\ldots,x_j)$ is {\em $(j,y)$-good} 
if $q_{j,y}(x_1,\ldots,x_j) \le 2q(j,y)$, and note
that for an i.i.d.($\mu$) sequence $(\xi_i)$
\begin{eqnarray}
\Pr ( (\xi_1,\ldots,\xi_j) \mbox{ is not $(j,y)$-good })
&=& 
\Pr( q_{j,y} (\xi_1,\ldots,\xi_j) > 2q(j,y) ) \nonumber\\
&\le& 
\frac{  \Ex q_{j,y} (\xi_1,\ldots,\xi_j) }{2q(j,y)} \nonumber\\
&=& \frac{ q(j,y)}{2q(j,y)} = \frac{1}{2} . \label{12q}
\end{eqnarray}
The proof of the next lemma exploits a certain monotonicity property: the spacings between customers tend to be shorter than usual, 
for the customers who have just moved in a long wave.
\begin{Lemma}
\label{Lsm}
For  $\bx \in \mathbb{X}, x^\prime < x_{s-t_0+1}(t_0-1)     , 1 \le t_0 < s, \ j \ge 1$,
\[ \hspace*{-1.6in} 
\Pr \left[ (\xi_{s-t_0+1}(t_0),\ldots, \xi_{s-t_0+j}(t_0)) \mbox{ is $(j,y)$-good } \right.
\]
\[  \hspace*{0.6in} \left.
\vert \bX(t_0-1) = \bx, X_{s-t_0}(t_0) = x^\prime, W(t_0) > s-t_0+j \right]
\ge \frac{1}{2} . \]
\end{Lemma}

{\bf Proof.}
Given $ \bX(t_0-1) = \bx$, the event $\{ W(t_0) > s-t_0+j  \}$ is the event
$\{ \sum_{i=1}^k \xi_i(t_0) < x_{k+1} - c^+, 1 \le k \le s - t_0 + j\}$.
This implies that, given $ \bX(t_0-1) = \bx$ and $X_{s-t_0} (t_0)= x^\prime < x_{s-t+1}(t_0-1)$ 
(the latter is equivalent to $X^*_s(t_0) = x^\prime < X^*_s(t_0-1)$),
the event $\{ W(t_0) > s-t_0+j  \}$ is the event 
\begin{equation}
\{ x^\prime + \sum_{i=1}^k \xi_{s-t_0+i}(t_0)  < x_{k+1} - c^+, 1 \le k \le j\}.
\label{event-x}
\end{equation}
Write $(\hat{\xi}_i, 1 \le i \le j)$ for random variables with the conditional distribution 
of $(\xi_{s-t_0+i}(t_0), 1 \le i \le j)$ given event (\ref{event-x}). 
Write $(\xi_i, 1 \le i \le j)$ for an independent i.i.d.($\mu$) sequence.
The conditional probability of event (\ref{event-x}) given $\xi_{s-t_0+1}(t_0)$
is clearly a decreasing function of  $\xi_{s-t_0+1}(t_0)$.  
From this and Bayes rule we see that the distribution of $\hat{\xi}_1$ is stochastically smaller than $\mu$, 
in other words we can couple $\hat{\xi}_1$ and $\xi_1$ so that $\hat{\xi}_1 \le \xi_1$ a.s. 
Repeating this argument inductively on $i$, we can couple $(\hat{\xi}_i, 1 \le i \le j)$ 
and $(\xi_i, 1 \le i \le j)$ so that  $\hat{\xi}_i \le \xi_i, 1 \le i \le j$.
Now the function  $q_{j,y}(x_1,\ldots,x_j) $ is monotone increasing in each argument $x_i$, so 
\[ \mbox{ if $(\xi_i, 1 \le i \le j)$ is $(j,y)$-good then $(\hat{\xi}_i, 1 \le i \le j)$  is $(j,y)$-good.} \]
Now (\ref{12q}) implies that $(\hat{\xi}_i, 1 \le i \le j)$  is $(j,y)$-good with probability $\ge 1/2$.
But this is the assertion of the lemma.
\qed

The next lemma starts to address the ``key idea" stated at the start of section \ref{sec:UB}.  
Given that a time-$t_0$ wave extends at least $j$ places past individual $s$, 
and given that a time-$t_0+t$ wave reaches individual $s$, 
we can upper bound the probability this time-$t_0+t$ wave extends at least $j$ places past individual $s$,
and the upper bound is small if individual $s$ has moved only distance $o(j^{1/2})$ between times $t_0$ 
and $t_0 +t$.
\begin{Lemma}
\label{Lt0}
For $t_0 \ge 1, t \ge 1, s > t_0+t,  j \ge 1$, 
\[  \hspace*{-0.9in}  \Pr(W(t_0+t) > s-t_0-t+j \vert \bX^*_{[s]}(t_0+t), \FF(t_0+t-1)) 
\] \[ 
\hspace*{0.9in} 
\le q_{j,X^*_s(t_0) - X^*_s(t_0+t)} (\xi_{s-t_0+1}(t_0),\ldots, \xi_{s-t_0+j}(t_0))
\] \[  \hspace*{0.9in} 
\mbox{ on } \{ W(t_0+t) > s-t_0-t \}  \cap \{ W(t_0) > s-t_0+j\}. \]
\end{Lemma}
Note that the event $\{ W(t_0+t) > s-t_0-t \} $ is the event 
$\{X^*_s(t_0+t) < X^*_s(t_0+t-1)\}$
and so is indeed in $\sigma (\bX^*_{[s]}(t_0+t), \FF(t_0+t-1))$. 
Also on the event $ \{ W(t_0) > s-t_0+j\}$ we have (for $1 \le i \le j$) that 
$\xi_{s-t_0+i}(t_0) = X_{s-t_0+i}(t_0) - X_{s-t_0+i-1}(t_0)$ and so 
$\xi_{s-t_0+i}(t_0)$ is $\FF(t_0)$-measurable.

{\bf Proof.}
On the event $\{ W(t_0+t) > s-t_0-t \} $ we have $\{X^*_s(t_0+t) < X^*_s(t_0+t-1)\}$, and given this, the 
event $\{W(t_0+t) > s-t_0-t+j \}$ is the event
\[
\{ X^*_s(t_0+t) + \sum_{i=1}^k \xi_{s-t_0-t+i}(t_0+t)  < X^*_{s+k}(t_0+t-1) - c^+, 1 \le k \le j\}.
\]
Now $X^*_{s+k}(t_0+t-1) \le X^*_{s+k}(t_0) $ and $-c^+ < 0$, so the event above is 
a subset of the event
\begin{equation}
\{ \sum_{i=1}^k \xi_{s-t_0-t+i}(t_0+t)  < X^*_{s+k}(t_0) -   X^*_s(t_0+t),   1 \le k \le j\} .
\label{iwt}
\end{equation}
Restricting further to the event $ \{ W(t_0) > s-t_0+j\}$, 
on which $X^*_{s+k}(t_0) - X^*_s(t_0) = \sum_{i=1}^k \xi_{s-t_0+i}(t_0), \ 1 \le k \le j$,
we can rewrite (\ref{iwt}) as the event
\begin{equation}
 \{ \max_{1 \le k \le j}  \sum_{i=1}^k  (\xi_{s-t_0-t+i}(t_0+t) - \xi_{s-t_0+i}(t_0)) < X^*_{s}(t_0) -   X^*_s(t_0+t) \} . 
 \label{event-x2}
 \end{equation}
So the conditional probability in the statement of the lemma is bounded by the conditional probability, 
given $\HH := \sigma( \bX^*_{[s]}(t_0+t), \FF(t_0+t-1))$, of  event (\ref{event-x2}).
The random variables $ \xi_{s-t_0+i}(t_0)), \ X^*_{s}(t_0), \ X^*_s(t_0+t) $ are $\HH$-measurable and the 
$\xi_{s-t_0-t+i}(t_0+t)$ are independent of them.
So the conditional probability of  event (\ref{event-x2}) given $\HH$ on the subsets specified equals precisely 
the bound stated in the lemma, by definition (\ref{def-qjyx}) of $q_{j,y}(\cdot)$.
\qed

The way we will combine these lemmas is abstracted in the next  lemma, in which
$S \wedge T$ denotes $\min(S,T)$.
Also, we rather pedantically write $\sigma(A_t, \FF(t-1))$ in the conditioning statement to emphasize that 
we are conditioning on a $\sigma$-field rather than an event.

\begin{Lemma}
\label{LAB}
Let $1 \le T < \infty$ be a stopping time for a filtration $(\FF(t), t = 0, 1,2,\ldots )$.  
Let $A_t$ and $B_t$ be events in $\FF_t$ satisfying
\begin{eqnarray}
B_t \subseteq A_t && \label{BAt} \\
\Pr(B_t \vert \sigma(A_t, \FF(t-1))) &\le& \delta \mbox{ on } A_t \cap \{T \ge t\}, \ t \ge 1. \label{BAF}
\end{eqnarray}
Define 
\[ N(\tau) = \sum_{t=1}^\tau \ind_{A_t} \]
\[ S = \min \{t \ge1 : \ B_t \ \mathrm{occurs }\}. \]
Then $\Pr(T \le S) \le \delta \Ex N(S \wedge T)$.
\end{Lemma}
{\bf Proof.}
\begin{eqnarray}
\Pr(B_t \vert \FF(t-1)) \ind_{\{T\ge t\}} &=&
\Ex \left[ \Pr(B_t \vert \sigma(A_t, \FF(t-1))) \ind_{\{T\ge t\}}  \ \vert \ \FF(t-1)\right] \nonumber \\
&\le& \Ex \left[ \delta \Pr(A_t \vert \FF(t-1)) \ind_{\{T\ge t\} } \ \vert \ \FF(t-1) \right] \nonumber \\
&=& \delta  \Pr(A_t \vert \FF(t-1)) \ind_{\{T\ge t\}}  \label{BdA}
\end{eqnarray}
the inequality holding  by assumptions (\ref{BAt}, \ref{BAF}).
Now 
\begin{eqnarray*}
\Pr(S = t, S \le T)
&=& \Pr(B_t, S \ge t, T \ge t) \\
&=& \Ex \left[  \Pr(B_t \vert \FF(t-1)) \ind_{\{T\ge t\}} \ind_{\{S\ge t\}}  \right] \\
&\le&  \delta \Ex \left[  \Pr(A_t \vert \FF(t-1)) \ind_{\{T\ge t\}}  \ind_{\{S\ge t\}} \right] \mbox{ by } (\ref{BdA})\\
&=& \delta \Pr(A_t, S \wedge T \ge t) .
\end{eqnarray*}
Summing over $t \ge 1$ gives the result.
\qed

\bigskip
We will now fit together the results above to obtain a result stated as Proposition \ref{P:Xst}  below. 
The statement involves substantial notation, and it is less cumbersome to develop the notation as we proceed.

Fix even $j \ge 4$ and $s>j$ and $y>0$. 
Start the process with  $\bX(0) = \bx$ arbitrary.
Consider events, for $1 \le t < s$,

\medskip \noindent
 $A_t = \{W(t) > s-t\}$ \\
$C_t = \{W(t) > s-t+j\}$ \\
$B_t = C_t \cap \{ (\xi_{s-t+1}(t), \ldots, \xi_{s-t+j}) \mbox{  is $(j,y)$-good } \}$.

In words:\\
$A_t = $ ``customer $s$ moves at time $t$" \\
$C_t = $ ``a long wave at time $t$"\\
 $B_t = $ ``the long wave at time $t$ is good".

\smallskip \noindent
Note that we are measuring wave length relative to the current position of customer $s$.
We  study what happens between the times of good long waves, that is the times
\begin{eqnarray}
T_1 &=& \min\{t \ge 1: \ B_t \mbox{ occurs } \} \\
T_u &=& \min\{t > T_{u-1}: \ B_t \mbox{ occurs } \}, u \geq 2 .
\end{eqnarray}
Note $T_u$ may be infinite (if no such event occurs); 
this will not be an issue because we will be conditioning on $\{T_u = t_0\}$ and will truncate $T_{u+1}$.

Fix $u \ge 1$ and $1 \le t_0 < j/2$. 
Recall we previously fixed $s > j \ge 4$ and $y>0$; all these variables are integers except $y$.

By definition,  $B_{t_0}$ is the event that the realization of
$(\xi_{s-t_0+i}(t_0) = X_{s-t_0+i}(t_0) - X_{s-t_0+i-1}(t_0), \ 1 \le i \le j)$
is $(j,y)$-good.
For each $1 \le t < s-t_0$, Lemma \ref{Lt0} and the fact that 
$q_{j,y}(\cdot)$ is monotone increasing in $y$ imply
\[ \hspace*{-0.8in}
   \Pr(C_{t_0+t} \vert \bX^*_{[s]}(t_0+t), \FF(t_0+t-1)) 
\leq 2q(j,y) \]
\begin{equation}
\mbox{ on }
A_{t_0+t} \cap B_{t_0}
\cap \{X^*_s(t_0) - X^*_s(t_0+t) \le y\}  \cap \{T_u = t_0\}. \label{Ctt}
\end{equation}
Now write 
\[ \widehat{T} = \min\{t \ge 1: X^*_s(t_0) - X^*_s(t_0+t) > y\} \]
so that the penultimate event in the intersections in (\ref{Ctt}) is the event $\{ \widehat{T} > t_0+t\}$.
Because $B_{t_0+t}  \subseteq C_{t_0+t} $ we now have
\[   \Pr(B_{t_0+t} \vert \bX^*_{[s]}(t_0+t), \FF(t_0+t-1)) 
\leq 2q(j,y) \mbox{ on }
A_{t_0+t} \cap B_{t_0} \cap \{ \widehat{T} > t_0+t\}  \cap \{T_u = t_0\}.
\]
Because $B_{t_0} \cap \{ \widehat{T} > t_0+t\}  \in  \FF(t_0+t-1)$  and $A_{t_0+t}  \in \sigma( \bX^*_{[s]}(t_0+t), \FF(t_0+t-1))$
we  can condition down from $\sigma( \bX^*_{[s]}(t_0+t), \FF(t_0+t-1))$ 
to $\sigma( A_{t_0+t}, \FF(t_0+t-1))$ and write
\begin{equation}
  \Pr(B_{t_0+t} \vert A_{t_0+t}, \FF(t_0+t-1)) 
\leq 2q(j,y) \mbox{ on }
A_{t_0+t} \cap B_{t_0} \cap \{ \widehat{T} > t_0+t\}  \cap \{T_u = t_0\}.
\label{BtAt}
\end{equation}
Define 
\[N_A(\tau) = \sum_{t=1}^\tau \ind_{A_t}\]
 and similarly for $N_B(\tau)$ and $N_C(\tau)$.

We can now apply Lemma \ref{LAB}, restated for times $t_0 < t < \infty$.
Condition on $\sigma( \{T_u = t_0\}, \FF(t_0-1), X^*_s(t_0))$ and restrict to the event $\{T_u = t_0\}$. 
In the hypotheses of  Lemma \ref{LAB} take   $T = \widehat{T}\wedge j/2, S = T_{u+1}, \delta =  2q(j,y)$.  
So inequality (\ref{BAF}) holds by (\ref{BtAt}).
The conclusion of that lemma 
is now
\begin{quote}
conditional on 
$\sigma( \{T_u = t_0\}, \FF(t_0-1), X^*_s(t_0))$ and restricted to the event $\{T_u = t_0\}$, we have;

 $\Pr(\widehat{T}\wedge j/2  \le T_{u+1}) \le  2q(j,y)  \Ex(N_A(T_{u+1} \wedge \widehat{T} \wedge j/2 ) - N_A(T_u))$.
\end{quote}
Rewriting this with all conditioning made explicit:
\[ \hspace*{-0.8in} 
\Pr(\widehat{T}\wedge j/2   \le T_{u+1} \vert  \{T_u = t_0\}, \FF(t_0-1), X^*_s(t_0)) 
\] \[
\le  2q(j,y)  \Ex(N_A(T_{u+1} \wedge \widehat{T}\wedge j/2 ) - N_A(T_u) \vert   \{T_u = t_0\}, \FF(t_0-1), X^*_s(t_0)  ) 
\mbox{ on } \{T_u = t_0\} . \]
The event $\{ \widehat{T} \wedge j/2  \le T_{u+1}\}$ is the event 
$\{T_{u+1} \le j/2, X^*_s(T_u) - X^*_s(T_{u+1}) \le y\}
\cup \{T_{u+1} > j/2\}$, 
and because $N_A(T_{u+1} \wedge \hat{T} \wedge j/2 ) \le N_A(T_{u+1} \wedge j/2)$ we have
\[ \Pr(  X^*_s(T_u) - X^*_s(T_{u+1}) \le y, T_{u+1} \le j/2   \vert  \{T_u = t_0\}, \FF(t_0-1), X^*_s(t_0)) 
\] \[
\le  2q(j,y)  \Ex(N_A(T_{u+1} \wedge j/2 ) - N_A(T_u) \vert   \{T_u = t_0\}, \FF(t_0-1), X^*_s(t_0)  ) 
\mbox{ on } \{T_u = t_0\} . \]
The notation has become rather cumbersome; we will make it more concise
by defining the $\sigma$-field $\FF^-(T_u)$ as follows.
\begin{quote}
For each $t_0$, the restriction of  $\FF^-(T_u)$ to $\{T_u = t_0\}$ is 
 the restriction of  $\sigma(  \{T_u = t_0\}, \FF(t_0-1), X^*_s(t_0) )$ to $\{T_u = t_0\}$.
\end{quote} 
Here $\FF^-(T_u)$ is a sub-$\sigma$-field of the usual pre-$T_u$ $\sigma$-field $\FF(T_u)$.
In words, when $T_u = t_0$ the usual $\FF(T_u)$ tells us all about $\bX(t_0)$ whereas $\FF^-(T_u)$
 only tells us $X^*_s(t_0)$ and the fact that a good long wave has just occured.
Note that $T_u$ is indeed $\FF^-(T_u)$-measurable, as is $X^*_s(T_u)$.
Now the inequality above can be rewritten more concisely as
\begin{Proposition}
\label{P:Xst} 
For $u \ge 1$, even $j \ge 4, s>j, y > 0$
\[  \Pr(  X^*_s(T_u) - X^*_s(T_{u+1}) \le y, T_{u+1} \le j/2  \vert \FF^-(T_u)) 
\] \[
\le  2q(j,y)  \Ex(N_A(T_{u+1}\wedge j/2) - N_A(T_u) \vert   \FF^-(T_u)) 
\mbox{ on } \{T_u < j/2\} . \]
\end{Proposition}
In the next section we will use this result with $y = o(j^{1/2})$, implying $q(j,y)$ is small.
Note that Proposition \ref{P:Xst} formalizes the property stated at the start of this section
\begin{quote}
 after one long wave, and before individual $s$ 
moves a distance $o(j^{1/2})$:

\smallskip \noindent
either (i)  there is unlikely to be another long wave; \\
or (ii) there are likely to be many waves that move individual $s$.
\end{quote}
except that the Proposition refers to {\em good} long waves.  
However, Lemma \ref{Lsm} says that each long wave has probability $\ge 1/2$ to be good, 
independently of the past, implying that the mean 
number of long waves between successive good long waves is $\le 2$, and more precisely
\begin{Lemma}
\label{P:NCA} 
For even $j \ge 4, \ s > j$
\[ \Ex N_C(j/2) \le 2 \Ex N_B(j/2) . \]
\end{Lemma}

\subsection{The upper bound}
\label{sec:ub}

Write
\begin{equation}
N_j(\tau) = \sum_{t = 1}^\tau \ind_{\{W(t) > j\}} 
\label{Njt}
\end{equation}
for the number of waves of length $> j$ up to time $\tau$, and then write
\[ \rho^+(j) = \limsup_{\tau \to \infty} \tau^{-1} \Ex N_j(\tau) . \]

\begin{Proposition}
\label{P:ub}
$\rho^+(j) = O(j^{-1/2})$ as $j \to \infty$.
\end{Proposition}

{\bf Outline of proof.}
We track the initial rank-$j$ customer for time $j/2$. 
If the number of good long waves  during this time is $O(j^{1/2})$ then we are done, because by Lemma \ref{P:NCA} 
the number of waves of length $>2j$ must also be $O(j^{1/2})$ and hence their rate is $O(j^{-1/2})$.
So suppose the number $N^*_j$ of good long waves is larger than  $O(j^{1/2})$.  Then we can choose 
$y_j = o(j^{1/2})$ such that $N^*_jy_j$ grows faster than $O(j)$.
Then for most of these $N^*_j$ waves the probability of 
$\{ X^*_s(T_u) - X^*_s(T_{u+1}) \le y_j  \}$ must be $> 1/2$ 
and so by Proposition \ref{P:Xst}, because $q(j,y_j) = o(1)$, 
 the mean number of waves of length $> j/2$ is a large multiple of 
  the mean number of waves of length $> 2j$. 
  In other words, writing $W$ for typical wave length,
  \[ \mbox{ if $\Pr(W > j) \neq O(j^{-1/2})$ then $\Pr(W>2j)/\Pr(W > j/2)$ is small } \]
 and this easily implies that in fact $\Pr(W > j) = O(j^{-1/2})$.

{\bf Details of proof.}
As in the previous section we track customer $s$, but we now set $s = j$ and track the customer for time $j/2$. 
Fix $y$. 
Take expectation in Proposition \ref{P:Xst} over $\{T_u < j/2\}$ to get
\[ \Pr(X^*_s(T_u) - X^*_s(T_{u+1}) \le y, \ T_{u+1}\le j/2) \]
\begin{equation}
\le 2 q(j,y) \Ex (N_A(T_{u+1}\wedge j/2) - N_A(T_u) ) \ind_{(T_u < j/2)} .
\label{qNa}
\end{equation}
Note that the number $U_j$ of $(j,y)$-good long waves before time $j/2$ can be represented as
\[ U_j = \min \{u: T_{u+1} > j/2\} . \] 
Summing (\ref{qNa}) over $u \ge 1$ gives
\begin{equation} 
\Ex \ | \{1 \le u \le U_j -1 : \ X^*_s(T_u) - X^*_s(T_{u+1}) \le y_j\}| 
\le 2 q(j,y) \   \Ex(N_A(j/2) - N_A(T_1) ).
\label{uUj}
\end{equation}
Because the initial position of customer $j$ is at most $c^+ j$ we have
\[ \sum_{u=1}^{U_j-1} ( X^*_s(T_u) - X^*_s(T_{u+1})) \le c^+ j \]
and hence
\[ | \{1 \le u \le U_j -1: \ X^*_s(T_u) - X^*_s(T_{u+1}) > y\}| \le c^+ j/y\]
implying
\[ | \{1 \le u \le U_j -1 : \ X^*_s(T_u) - X^*_s(T_{u+1}) \le y \}| \ge U_j -1 - c^+ j/y .\]
Taking expectation and combining with (\ref{uUj}), and noting $U_j = N_B(j/2)$, we see
\[ \Ex N_B(j/2) - 1 - c^+ j/y \le 2 q(j,y) \   \Ex N_A(j/2) . \]
Now our customer $s$ has rank decreasing from $j$ to $j/2$, so a wave that moves this customer must have length at least $j/2$:
\[ N_A(j/2) \le N_{j/2}(j/2) \]
in the notation at (\ref{Njt}). 
Similarly, a wave reaching rank $2j$ is a long wave, so $N_C(j/2) \ge N_{2j}(j/2)$, 
and combining with Lemma \ref{P:NCA} shows
\[ \Ex N_B(j/2) \ge \sfrac{1}{2} \Ex N_{2j}(j/2) .\]
So now we have
\[ \Ex N_{2j}(j/2) \le 2 + 2c^+ j/y  + 4 q(j,y) \Ex N_{j/2}(j/2) .\]
The initial configuration is arbitrary, so by averaging over consecutive time-intervals of length $j/2$ and
 taking time-asymptotics, 
  \[ \sfrac{1}{4} \rho^+(2j) \le \sfrac{1}{j} + \sfrac{c^+}{y} + q(j,y)  \rho^+(j/2) . \]
 Now set
 \[ \psi(j) = j^{1/2} \rho^+(j) . \]
Using the Lemma \ref{LRW} bound $q(j,y) \le Cy j^{-1/2} + \delta_j$, 
 the previous inequality becomes
 \[ \sfrac{1}{4} \psi(2j) \le \sfrac{2^{1/2}}{j^{1/2}} + \sfrac{2^{1/2}j^{1/2}c^+}{y} 
 + \sfrac{2Cy}{j^{1/2}} \psi(j/2) + 2 \delta_j \psi(j/2) . \]
 Optimizing over choice of $y = y_j$ gives
 \[ 
\sfrac{1}{4} \psi(2j) \le \sfrac{2^{1/2}}{j^{1/2}} + 2^{7/4}(c^+C)^{1/2} \psi^{1/2}(j/2) + 2 \delta_j \psi(j/2)  . \]
From this it is clear that
\[ \limsup_{i \to \infty} \psi(4^i) < \infty \]
and then because $\rho^+(j)$ is decreasing we have $\limsup_{i \to \infty} \psi(j) < \infty$.
\qed

\newpage
\section{Lower bounding the rate of long waves}
\label{sec:LB}

\subsection{The blocks argument}
Here we will employ the ``blocks" argument mentioned in the introduction. 
Fix $t_0 > 1$ and $i \ge 0$ and consider individual $s = t_0+i$, who is at position
$X_{i}(t_0) = X^*_s(t_0)$ at time $t_0$.  
Looking backwards in time, the time since the last move of that individual is
\[
L_s(t_0) = \max\{ \ell  \ge 0: X^*_s(t_0) = X^*_s(t_0 - \ell)\} 
\]
where we set  $L_s(t_0) = t_0$ if $s$ has never moved, that is if $X^*_s(t_0) = X^*_s(0)$.
Recall the notation
 \[ \bX_{[s]}(t) = (X_i(t), \ 0 \le i \le s-t) \]
 for the positions of individual $s$ and the customers ahead of that individual in the queue.
Consider the event 
\begin{equation}
A_{s+1}(t_0) := \{L_s(t_0) = L_{s+1}(t_0) < t_0\}
\label{def-As}
\end{equation}
that individuals $s$ and $s+1$ are in the same ``block" of customers whose last move before $t_0$ was at the same time.
Define
\[ 
\GG_s(t_0) = \sigma( \bX_{[s]}(t), 0 \le t \le t_0; A_{s+1}(t_0) ) .
\]

\begin{Lemma}
\label{L:almost}
For $t_0 > 1$ and $i \ge 0$ and $s = t_0+i$,
\[ 
\Pr(X^*_{s+1}(t_0) - X^*_{s} (t_0) \in \cdot \vert \GG_s(t_0)) = \mu(\cdot) \mbox{ on } A_{s+1} (t_0)
. \]
\end{Lemma}
So by conditioning down and changing notation to $X_i(t_0) = X^*_s(t_0)$  we have
\begin{equation}
\Pr(X_{i+1}(t_0) - X_{i} (t_0) \in \cdot \vert  X_j(t_0), 1 \le j \le i,   A_{t_0+i+1}(t_0)  ) = \mu(\cdot) \mbox{ on } A_{t_0+i+1} (t_0) .
\label{XitX}
\end{equation}
We will later show that  $1 - \Pr(A_{t_0+i+1} (t_0)) $ becomes small for large $i$, and then
 Lemma \ref{L:almost} provides one formalization of the idea, mentioned earlier, 
that the  inter-customer distances within the queue  are roughly i.i.d. ($\mu$).

{\bf Proof  of Lemma \ref{L:almost}.}
The history of the individual $s$, that is 
$(X^*_s(t), 0 \le t \le t_0)$, tells us the value $\ell = L_s(t_0)$.
Consider that history and the history of 
individual $s+1$ up to time $t_0-\ell -1$.
Individual $s+1$ moved at time $t_0 - \ell$ 
if and only if the event  
\[ B_\ell := \{X^*_s(t_0 - \ell) < X^*_{s+1}(t_0 - \ell - 1) - c^+ \} \]
 happened. 
On the event $B_\ell$, by construction (\ref{Xit1} -  \ref{Xit3}) 
$X^*_{s+1}(t_0 - \ell) - X^*_s(t_0-\ell) = \xi_{i+\ell+1}(t_0 - \ell)$.
But time $t_0-\ell$ was the time of the last move of  individual $s$, so individual $s+1$ cannot subsequently move; 
the spacing is unchanged.
See Figure 5.

\setlength{\unitlength}{0.03in}
\begin{picture}(100,97)(-20,-10)
\put(20,70){\circle*{1.91}}
\put(20,60){\circle*{1.91}}
\put(20,50){\circle*{1.91}}
\put(20,40){\circle*{1.91}}
\put(30,30){\circle*{1.91}}
\put(30,20){\circle*{1.91}}
\put(30,10){\circle*{1.91}}
\put(40,0){\circle*{1.91}}
\put(37,30){\circle*{1.91}}
\put(37,20){\circle*{1.91}}
\put(37,10){\circle*{1.91}}
\put(52,0){\circle*{1.91}}
\put(20,70){\line(0,-1){30}}
\put(20,40){\line(1,-1){10}}
\put(30,30){\line(0,-1){20}}
\put(30,10){\line(1,-1){10}}
\put(37,30){\line(0,-1){20}}
\put(37,10){\line(3,-2){15}}
\put(20,40){\vector(1,0){12}}
\put(24,42){$\xi =  \xi_{i+\ell+1}(t_0 - \ell)$}
\put(-10,69){$t_0$}
\put(-10,39){$t_0 - \ell$}
\put(-10,29){$t_0 - \ell -1$}
\put(-10,-2){$0$}
\put(-30,54){time}
\put(5,-5){individuals} 
\put(38,-5){$s$}
\put(49,-5){$s +1$}
\put(-15,0){\line(0,1){80}}
\put(-15,0){\line(1,0){2}}
\put(-15,30){\line(1,0){2}}
\put(-15,40){\line(1,0){2}}
\put(-15,70){\line(1,0){2}}
\end{picture}

{\bf Figure 5.} Individual $s+1$ moves at time $t_0 - \ell$ on event $B_\ell$.

\medskip \noindent
So
\[ X^*_{s+1}(t_0) - X^*_{s} (t_0) =  \xi_{i+\ell+1}(t_0 - \ell) \mbox{ on } B_\ell \cap  \{L_s(t_0) = \ell\} . \] 
Now the event $B_\ell \cap  \{L_s(t_0) = \ell\} $ 
and the random variables 
$(\bX_{[s]}(t), 0 \le t \le t_0, 
X_{s+1}(t), 0 \le t \le t_0-\ell - 1) $
are determined by members of the constructing family $(\xi_j(t), j,t \ge 1)$ 
not including $ \xi_{i+\ell+1}(t_0 - \ell)$.
So by independence of the constructing family,
\[ \Pr(X^*_{s+1}(t_0) - X^*_{s} (t_0) \in \cdot \vert \bX_{[s]}(t), 0 \le t \le t_0, 
X_{s+1}(t), 0 \le t \le t_0-\ell - 1) 
\] \[  = \mu(\cdot) \mbox{ on } B_\ell \cap  \{L_s(t_0) = \ell\} . 
\]
Now 
$A_{s+1}(t_0)$ is the disjoint union $ \cup_{0 \le \ell < t_0}  [  B_\ell \cap  \{L_s(t_0) = \ell\}   ] $
and the result follows easily.
\qed

{\bf Remark.} 
Using Lemma \ref{L:almost} and (\ref{XitX}) to deduce properties of the distribution of 
$\bX(t_0) = (X_i(t_0), i \ge 0)$ requires some care, because
there is some complicated dependence between the event $A_{t_0+i+1}(t_0)$ and the tail sequence 
$(X_j(t_0), j \ge i+1)$ which we are unable to analyze.
This is the obstacle to a
``coupling from the past" proof that for large $t_0$ the configuration $\bX(t_0)$ can be defined in terms of the 
past $\xi$'s with vanishing dependence on the initial configuration $\bX(0)$, and thereby proving weak convergence to a unique stationary distribution.
 
The following argument is unaffected by such dependence.
Fix $t_0$ and inductively on $i \ge 0$ construct $\xi^*_{i+1}(t_0)$ by

\noindent
$ \mbox{ on } A_{t_0+i+1}(t_0)$ let $ \xi^*_{i+1}(t_0)= X_{i+1}(t_0) - X_{i} (t_0) $

\noindent
$ \mbox{ on } A^c_{t_0+i+1}(t_0) $ take $\xi^*_{i+1}(t_0)$ to have distribution $\mu$ independent of all previously-defined random variables.

Then using (\ref{XitX}) the sequence $(\xi^*_i(t_0), i \ge 1)$ is i.i.d.($\mu$).
Write $S_k(t_0) = \sum_{i=1}^k \xi^*_i(t_0)$.
By construction we immediately have
\begin{Lemma}
\label{L:XkS}
\[ |X_k(t_0) - S_k(t_0)| \le (c^+ - c_-) \sum_{i=1}^k \ind_{A^c_{t_0+i}(t_0)} . \]
\end{Lemma}

\subsection{The lower bound}
\label{sec:lb}
We can now outline the proof that waves of length $> k$ occur at rate no less than order $k^{-1/2}$. 
By showing that the bound in Lemma \ref{L:XkS} is sufficiently small (Lemma \ref{L:XkU}), 
we will see that the rank-$k$ individual at time $t_0$, for typical large $k$ and $t_0$, is at position 
$k \pm O(k^{1/2})$. 
By the considering the same estimate for the same individual a time $D k^{1/2}$ later (for large $D$), 
the individual must likely have moved during that time interval.  
In other words, there is likely to have been a wave of length $> k$ during this $O(k^{1/2})$ 
time interval.

In section \ref{sec:PP7} we will prove the bound in the following form.
Let $U_\tau$ have uniform distribution on $\{1,2,\ldots,\tau\}$,
\begin{Proposition}
\label{P:lb}
\[ \limsup_{\tau \to \infty} 
\Pr \left(\max_{U_\tau  < t \le U_{\tau }+ j  } W(s) \le k \right) \le Bk^{1/2}/j , \quad j < k  \]
for constant $B$ not depending on $j,k$.
\end{Proposition}
To translate this into the same format as Proposition \ref{P:ub}, write
\[ \rho_-(j) = \liminf_{\tau \to \infty} \tau^{-1} \Ex N_j(\tau) \]
where as before
\[
N_j(\tau) = \sum_{t = 1}^\tau \ind_{\{W(t) > j\}} 
\]
is the number of waves of length $> j$ up to time $\tau$.
\begin{Corollary}
\label{C:lb}
$\rho_-(k) \ge \frac{1}{2(1 + Bk^{1/2})}$, for the constant $B$  in Proposition \ref{P:lb}.
\end{Corollary}
{\bf Proof of Corollary \ref{C:lb}.}
The event in Proposition \ref{P:lb} is the event 
$\{ N_k(U_\tau + j) = N_k(U_\tau)\}$, so choosing 
$j = \lceil 2 Bk^{1/2} \rceil$ we have
\[ \liminf_{\tau \to \infty} 
\Pr (N_k(U_\tau + j) - N_k(U_\tau) \ge 1) \ge \sfrac{1}{2} \]
implying
\[ \liminf_{\tau \to \infty} 
\Ex(N_k(U_\tau + j) - N_k(U_\tau) ) \ge \sfrac{1}{2}. \]
For $1 \le t \le \tau+j$ we have 
$\Pr(U_\tau < t  \le U_\tau + j) \le j/\tau$; 
applying this to the times $t$ counted by $N_k(\tau + j)$ gives
\[ \Ex(N_k(U_\tau + j) - N_k(U_\tau) ) \le \sfrac{j}{\tau} \Ex N_k(\tau +j) . \]
So 
\[ \liminf_{\tau \to \infty} \sfrac{1}{\tau} \Ex N_k(\tau +j) \ge \sfrac{1}{2j} \]
establishing the Corollary.
\qed

\subsection{Proof of Proposition \ref{P:lb}}
\label{sec:PP7}
In order to apply Lemma \ref{L:XkS} we need to upper bound the probability of the complement of events 
$A_{s+1}(t_0)$ of the form (\ref{def-As}).
From the definition
\[ A^c_{s+1}(t_0) = \{ L_s(t_0) = \ell < L_{s+1}(t_0) \mbox{ for some } 0 \le \ell < t_0\}
\cup \{L_s(t_0) = t_0\} . \]
In order that event $\{ L_s(t_0) = \ell < L_{s+1}(t_0) \}$ occurs it is necessary 
(but not sufficient -- see remark below Lemma \ref{L:XkU}) that the wave at time $t_0 - \ell$ moves individual $s$ but not 
individual $s+1$, which is saying that $W(t_0 - \ell) = s - t_0 + \ell +1$.
And the event $\{ L_s(t_0) = t_0\} $ can be rewritten as  the event 
$\{M(t_0) \le s \}$ for
\[
M(t_0) := \min \{s: X^*_s(t_0) = X^*_s(0)\} .
\]
So
\[
\Pr ( A^c_{s+1}(t_0) ) \le \Pr (  M(t_0) \le s ) 
+ \sum_{\ell = 0}^{t_0-1} \Pr (   W(t_0 - \ell) = s - t_0 + \ell +1 ) .
\]
Setting $s = t_0 + i - 1$ for $i \ge 1$ and 
$\ell = t_0 - j$, this becomes
\[ \Pr ( A^c_{t_0 +i}(t_0) ) \le \Pr (  M(t_0) \le t_0 +i-1 )  
+ \sum_{j=1}^{t_0} \Pr (W(j) = t_0 + i - j) . \]
Now note that if individual $s$ is at the same position $x$ at times $0$ and $t_0$ we must have
\[   c_- s \le x \le c^+(s- t_0) \]
implying that $s \ge \sfrac{t_0c^+}{c^+ - c_-}$.  So there is a deterministic bound
$ M(t) \ge \sfrac{t_0c^+}{c^+ - c_-}$
which implies that, in the limit we will take at (\ref{tautau}), the contribution from the term 
$\Pr (  M(t_0) \le t_0 +i-1 )  $ is negligible, so we ignore that term in the next calculation.
We now average over $1 \le t_0 \le \tau$:
\begin{eqnarray*}
 \tau^{-1} \sum_{t_0 = 1}^\tau \Pr ( A^c_{t_0 +i}(t_0) ) &\le &
\tau^{-1} \sum_{t_0 = 1}^\tau \sum_{j=1}^{t_0} \Pr (W(j) = t_0 + i - j) \\
&=&  \tau^{-1} \sum_{j = 1}^\tau \Pr( i \le W(j) \le \tau + i - j) \\
&\le&  \tau^{-1} \sum_{j = 1}^\tau \Pr( W(j) \ge i)  
\end{eqnarray*}
and take limits
\begin{equation}
 \limsup_{\tau \to \infty} \tau^{-1} \sum_{t_0 = 1}^\tau \Pr ( A^c_{t_0 +i}(t_0) ) \le 
\limsup_{\tau \to \infty}  \tau^{-1} \sum_{j = 1}^\tau \Pr( W(j) \ge i)  . 
\label{tautau}
\end{equation}
Proposition \ref{P:ub} bounds the right side as $O(i^{-1/2})$. 
Now we fix $k$, sum over $1 \le i \le k$ and apply Lemma \ref{L:XkS} to deduce
\begin{Lemma}
\label{L:XkU}
There exists a constant $C$ such that, for each $k \ge 1$,
\begin{equation}
 \limsup_{\tau \to \infty} 
\Ex  |X_k(U_\tau) - S_k(U_\tau)| \le C k^{1/2} .
\label{tau4}
\end{equation}
\end{Lemma}
Recall that  $U_\tau$ has uniform distribution on $\{1,2,\ldots,\tau\}$, and that $S_k(t_0)$, and hence 
$S_k(U_\tau)$, is distributed as the sum of $k$ i.i.d.($\mu$) random variables.

{\bf Remark.}
Heuristics  suggest that in fact $ \Pr ( A^c_{t_0 +i}(t_0)) $ decreases as order 
$i^{-1} \log i$.  Our order $i^{-1/2}$ bound is crude because, in order for event $A_{s+1}(t_0)$ to occur, we need not only that
 $L_s(t_0) = \ell < L_{s+1}(t_0) \mbox{ for some } 0 \le \ell < t_0$, but also that no subsequent wave during time 
 $(t_0 - \ell, t_0]$ moves individual $s+1$, and we ignored the latter condition.  
 However a better estimate would not help our subsequent arguments, because we use (\ref{tau4}) to deduce an 
 $O(k^{1/2})$ bound on the spread of $X_k(U_\tau)$ from the same order bound on the spread of $S_k(U_\tau)$, 
 and this deduced bound would not be improved by a better bound on the difference.

\medskip \noindent
We can now continue with the argument outlined at the start of section \ref{sec:lb}.
We already fixed $k$, and now fix $j < k$.
Consider the individual of rank $k$ at time $t_0$, and therefore of rank $k+j$ at time $t_0-j$.
We have
\begin{equation}
\Pr(\max_{t_0-j  < t \le t_0} W(t) \le k) \le \Pr( X_k(t_0) = X_{k+j}(t_0-j)) 
\label{tjW}
\end{equation}
because if there was no wave of length $> k$ during the time interval $(t_0 -j,t_0]$ then
the individual under consideration does not move.  
Write
\begin {eqnarray*} 
X_{k+j}(t_0-j) - j - X_k(t_0) 
&=& S_{k+j}(t_0-j) - (k+j) \\
&-& ( S_k(t_0) - k) \\
&+& ( X_{k+j}(t_0-j) - S_{k+j}(t_0-j) ) \\
&-& ( X_k(t_0) - S_k(t_0) ) .
\end{eqnarray*}
Now $S_k(t_0)$ has mean $k$ and variance $k \sigma^2$, 
and $S_{k+j}(t_0-j)$ has mean $k+j$ and variance $(k+j)\sigma^2$, 
so by bounding the expectation of absolute value of each of the four terms in the sum above,
\[
\Ex \left| X_{k+j}(t_0-j) - j - X_k(t_0) \right| 
\le 
\] \[
(k^{1/2} + (k+j)^{1/2}) \sigma +
\Ex \left| X_{k+j}(t_0-j) - S_{k+j}(t_0-j) \right|
+ \Ex \left| X_k(t_0) - S_k(t_0) \right| . \]
The limit assertion in (\ref{tau4}) remains true if $U_\tau$ is replaced therein by $U_\tau + j$,
and so substituting  $U_\tau + j$ for $t_0$ in the inequality above and applying (\ref{tau4}) 
\[  \limsup_{\tau \to \infty} 
\Ex \left| X_{k+j}(U_\tau ) - j - X_k(U_\tau +j) \right| 
\le 2(\sigma + C) (k+j)^{1/2} . \]
Then by Markov's inequality,
\[  \limsup_{\tau \to \infty} \Pr (X_{k+j}(U_\tau ) = X_k(U_\tau +j) )
\le \frac{ 2(\sigma + C) (k+j)^{1/2}}{j} . \]
Now substitute $U_\tau + j$ for $t_0$ in (\ref{tjW}) and combine with the inequality above:

\[ \limsup_{\tau \to \infty} 
\Pr \left(\max_{U_\tau  < t  \le U_\tau +j} W(t) \le k \right) 
\le \frac{ 2(\sigma + C) (k+j)^{1/2}}{j} . \]
establishing Proposition \ref{P:lb}.

\newpage
\section{The CBM limit}
\label{sec:CBMlimit}
Here we develop the idea, outlined in section \ref{sec:CBM}, that the representation of the spatial queue via the process 
$G(t,x)$ in Figure 3,  suitably rescaled, converges in some sense to CBM. 
The notation needed to state  carefully the result, Theorem \ref{T:main}, is provided in  section \ref{sec:abstract}. 
Given the results and techniques from sections \ref{sec:UB} and  \ref{sec:LB},
the proof of Theorem \ref{T:main} (section \ref{sec:outpf}) is conceptually straightforward,
and we will not provide all the details of checking various intuitively clear assertions, for instance 
those involving the topology 
of our space $\Fspace$.  
A technical discussion of topologies in more general contexts of this kind can be found in \cite{Bweb}, but here we set out a 
bare minimum required to formulate the result.
Finally in section \ref{sec:s-t} we outline how the results concerning customers' 
space-time trajectories follow from the CBM limit.

\subsection{Notation and statement of result}
\label{sec:abstract}

\medskip \noindent $\bullet$ \hspace{0.1in} 
$C[1, \infty)$ is the space of continuous functions $f: [1,\infty) \to \Reals$, and
$d_0$ is a metrization of the usual topology of uniform convergence on compact intervals.

\medskip \noindent $\bullet$ \hspace{0.1in} 
 $\bdf = \{f \in \II(\bdf)\}$ denotes a ``coalescing family" of functions $f \in  C[1, \infty)$   with some index set $\II(\bdf)$, 
satisfying the conditions (i, ii) following.

(i)   $\init(\bdf) := \{f(1), \ f \in \II(\bdf)\}$ is a distinct locally finite set;

\noindent
so if $f \in  \II(\bdf)$ is such that $f(1)$ is not maximal in $\init(\bdf)$, then there is a ``next" $f^* \in \II(\bdf)$, the function such that 
$f^*(1)$ is the next largest in  $\init(\bdf)$.  For such a pair $(f,f^*)$

(ii) $f^*(t) \ge f(t) \ \forall t \ge 1$, and there is a ``coalescing time" $\gamma(f,f^*) < \infty$ such that $f^*(t) = f(t) \ \forall t \ge \gamma(f,f^*)$.

\noindent
Note this implies that each pair $(f_1,f_2)$ of functions in $\bdf$ will eventually coalesce at some time 
$\gamma(f_1,f_2)$.

\medskip \noindent $\bullet$ \hspace{0.1in} 
$\Fspace$ is the space of such families $\bdf$, equipped with the ``natural" topology in which 
$\bdf^{(n)} \to \bdf$ means the following:

for each finite interval $I \subset \Reals$ such that  $\init(\bdf)$ contains neither endpoint of $I$, 
for all large $n$ there is an ``onto" map 
$\iota^{(n)}$ from $\{f^{(n)} \in \bdf^{(n)}: f^{(n)}(1) \in I\}$ to 
$\{f \in \bdf: f(1) \in I\}$ such that, 
for  all choices of $f^{(n)}  \in \bdf^{(n)}$, we have that
$f^{(n)}    - \iota^{(n)} (f^{(n)}) $ converges to the zero function in $C[1, \infty)$; and also 
there is convergence of the coalescing times, in that
$\gamma^{(n)}(f^{(n)}_1, f^{(n)}_2) - \gamma(\iota^{(n)}(f^{(n)}_1), \ \iota^{(n)}(f^{(n)}_2)) \to 0$.

(Note that this definition is designed to allow there to be two functions in $\bdf^{(n)}$ converging to the 
same function in $\bdf$, provided their coalescence time $\to 1$.)


\medskip \noindent $\bullet$ \hspace{0.1in} 
For $\bdf \in \Fspace$ and $t_0 > 1$, write $\bdf_{[t_0]}$ for the time-shifted and rescaled family in which each
$f \in \bdf$ is replaced by $f^\prime$ defined by 
\begin{equation}
f^\prime(t) = t_0^{-1/2} f(t_0 t), \ t \ge 1. 
\label{def-fp}
\end{equation}

\medskip \noindent $\bullet$ \hspace{0.1in} 
Write $\bF$ for a random family, that is a random element of $\Fspace$, and 
$\bF_{[t_0]}$ for the time-shifted and rescaled family.

\medskip \noindent $\bullet$ \hspace{0.1in} 
Write $\CB$ for the specific random family consisting of the standard CBM process observed  from time $1$ onwards.
That is, the family of distinct functions
$ (t \to B(y,t) , t \ge 1)$ as $y$ varies.  
The scale-invariance property of standard CBM implies that $\CB_{[t_0]}$ has the same distribution as $\CB$.

\bigskip \noindent
Recall the discussion surrounding the Figure 3 graphic.
When this graphic represents the spatial queue, the horizontal axis is ``space" $x$ and the vertical axis is 
``time" $t$.  
To compare the spatial queue process with CBM we need to reconsider the horizontal axis as ``time" $t$ and the vertical axis as ``space" $y$.
That is, we rewrite the function $G(t,x)$ at (\ref{def-G}) as $\tilde{G}(y,t)$, so we have
\[
\tilde{G}(y,t) = G(y,t) = y + F_y(t) , \ \mbox{ where }
F_y(t) = \max\{k: X_k(y) \le t \} \]
where $y = 0,1,2,\ldots$ and $0 \le t < \infty$.
To study this process for large $y$ we rescale as follows:
for each $n$ define
\begin{equation}
\widetilde{H}^{(n)}(y,t) =  \frac{ \tilde{G}(\sigma n^{1/2}y + U_{\tau_n}, nt)  \ - \  U_{\tau_n} }{n^{1/2}}, 
 \quad 1 \le t < \infty .
\label{def-Ht}
\end{equation}
This is defined for $y$ such that $\sigma n^{1/2}y + U_{\tau_n}  \in \Ints^+$;
also $\tau_n \to \infty$ is specified below and $ U_{\tau_n} $ is uniform random on $\{1,\ldots,\tau_n\}$.
We now define the random family $\bF^{(n)}$ to consist of the distinct functions 
\[ t \to \widetilde{H}^{(n)}(y,t),    \quad 1 \le t < \infty  \]
as $y$ varies. 
This rescaling construction is illustrated in Figure 6.

\includegraphics[width=5.0in]{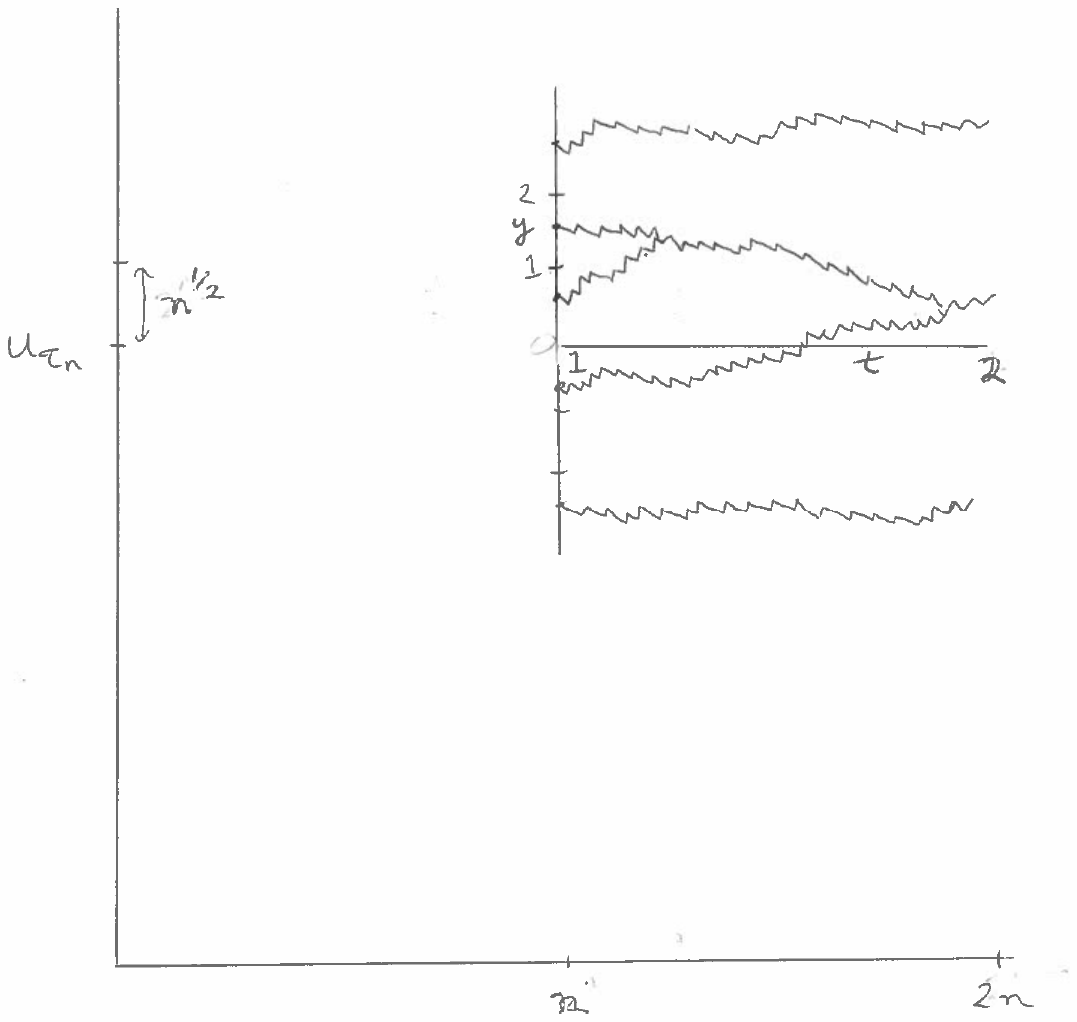}

\vspace*{-1.4in}

{\bf Figure 6.}  The rescaling that defines  $\bF^{(n)}$ (small axes) in terms of 
$G(t,x)$ (large axes), assuming $\sigma = 1$.

\bigskip

We can now state the convergence  theorem, in which
we are considering $\bF^{(n)}$ and $ \CB$ as random elements of the space $\Fspace$.
\begin{Theorem}
\label{T:main}
If $\tau_n \to \infty$ sufficiently fast, then $\bF^{(n)} \cd   \CB$.
\end{Theorem}
We conjecture this remains true without time-averaging, that is  if we replace $U_{\tau_n}$ by $\tau_n$ 
in the definition of $\bF^{(n)}$.

\subsection{Outline proof of Theorem \ref{T:main}}
\label{sec:outpf}

{\bf Step 1.} 
Propositions \ref{P:ub} and \ref{P:lb} refer to wave lengths in terms of rank $j$ instead of position (which we are now calling 
$t$).
But for any individual at any time these are related by $t/j \in [c_-,c^+]$, so the
``order of magnitude" bounds in those Propositions  remain true when we measure wave length by position.
By choosing $\tau_n \to \infty$ sufficiently fast, those Propositions provide information 
about the point processes 
$( \init(\bF^{(n)}) ,  n \ge 1)$, as follows.
\begin{Corollary}
There exist  constants $\beta_1, \beta_2$ and  $\zeta_n \to - \infty$ such that for all sufficiently large $n$
\begin{equation}
   \Pr ( \init(\bF^{(n)}) \cap [a, b] = \emptyset) 
\le \beta_1/(b-a) , \quad \zeta_n \le a  < b < \infty  .
   \label{Fnab1}
\end{equation}
\begin{equation}
   \Ex |  \init(\bF^{(n)}) \cap [a,b] | \le \beta_2 (b-a), \quad -\infty < a  < b < \infty  . 
   \label{Fnab2}
   \end{equation}
   \end{Corollary}

\bigskip \noindent {\bf Step 2.} 
Modify the usual topology on simple  point processes on $\Reals$ to allow two distinct points in the sequence to converge to the same point in the limit.  
In this topology, (\ref{Fnab2}) is a sufficient condition for tightness of the sequence 
$(\init(\bF^{(n)}), n \ge 1)$, and so by passing to a subsequence we may assume 
\begin{equation}
 \init(\bF^{(n)}) \cd \mbox{ some } \tilde{\eta}^{0} 
 \label{eta-P0}
 \end{equation}
where the limit point process $\tilde{\eta}^{0}$ inherits properties (\ref{Fnab1}, \ref{Fnab2}):
\[
   \Pr ( \tilde{\eta}^{0}  \cap [a, b] = \emptyset) 
\le \beta_1/(b-a) , \  - \infty < a  < b < \infty  .
\]
\begin{equation}
  \Ex |  \tilde{\eta}^{0}   \cap [a,b] | \le \beta_2 (b-a), \quad -\infty < a  < b < \infty  . 
  \label{eta-P2}
  \end{equation}
Now $\tilde{\eta}^{0} $ has translation-invariant distribution because we use the uniform random variable 
$U_{\tau_n}$ with $\tau_n \to \infty$ to define $\bF^{(n)}$.
The central part of the proof of Theorem \ref{T:main} is now to show
\begin{Proposition}
If $\tau_n \to \infty$ sufficiently fast then
\begin{equation}
 \bF^{(n)} \cd   \bF^{(\infty)} := \mbox{ CBM$(\tilde{\eta}^{0},0)$}
\label{FnC}
\end{equation}
where the limit is the coalescing Brownian motion process started with particle positions distributed as 
$\tilde{\eta}^{0}$.
\end{Proposition}
Granted (\ref{FnC}), we can fix $t_0 > 0$ and apply our
``time-shift and rescale" operation $\bdf \to \bdf_{[t_0]}$
 at (\ref{def-fp}) to deduce
 \begin{equation}
 \bF^{(n)}_{[t_0]} \cd   \bF^{(\infty)}_{[t_0]} \mbox{ as } n \to \infty . 
\label{FnC2}
\end{equation}
Now $ \bF^{(n)}$ depends on $\tau_n$, so to make that dependence explicit let us rewrite (\ref{FnC2}) as
\[  \bF^{(n, \tau_n)}_{[t_0]} \cd   \bF^{(\infty)}_{[t_0]} \mbox{ as } n \to \infty . \]
But $ \bF^{(n, \tau_n)}_{[t_0]} $ is just $\bF^{(nt_0)}$ with a different shift; precisely
\[ \bF^{(n, \tau_n)}_{[t_0]}  = \bF^{(nt_0, \tau_{nt_0})} . \] 
So we have 
\[ \bF^{(nt_0, \tau_{nt_0})} \cd  \bF^{(\infty)}_{[t_0]} \mbox{ as } n \to \infty . \]
But the scaling and embedding properties of CBM
imply that 
\[  \mbox{  $\bF^{(\infty)}_{[t_0]} \cd \CB$ as $t_0 \to \infty$ } . \] 
So by taking $t_0 = t_0(n) \to \infty$ sufficiently slowly and setting $m = n t_0(n)$,
\[ \bF^{(m,.\tau^\prime_m)} \cd \CB \mbox{ as } m \to \infty \]
for some $\tau^\prime_m \to \infty$, and this remains true for any larger $\tau^{\prime \prime}_m$.
This is the assertion of Theorem \ref{T:main}.

\bigskip \noindent {\bf Step 3.}
To start the proof of (\ref{FnC}), fix $Y_0$ and $T_0$, and 
define a modified process $\rmod_{Y_0}(\bF^{(\infty)})$ as CBM started with particles at positions 
$\tilde{\eta}^{0} \cap  [-Y_0,Y_0]$ only.
Define $\rmod_{Y_0}(\bF^{(n)})$ as follows.  
Consider the construction (\ref{Xit1} -  \ref{Xit3}) of the spatial queue process.
There is a smallest $t$ (in the notation of (\ref{Xit1} -  \ref{Xit3})) 
such that the corresponding value $y$ in $\init(\bF^{(n)})$, that is 
\[\sigma n^{1/2}(y + U_{\tau_n}) = t ,\]
satisfies $y \ge -Y_0$.  
We now modify  (in Step 4 we will argue that this modification has no effect in our asymptotic regime) 
the construction (\ref{Xit1} -  \ref{Xit3}) of the spatial queue process by saying that, for this particular $t$, we 
set
\[ X_i(t) =   \xi_1(t) + \ldots + \xi_i(t) , \  1 \le i < \infty . \]
That is, we replace all the inter-customer distances to the right of position $n$ by an i.i.d. sequence.
For subsequent times $t+1,t+2, \ldots$ 
we continue the inductive construction (\ref{Xit1} -  \ref{Xit3}). 
The effect of this change on the behavior of the queue process
 to the left of position $n$, in particular the effect on the process $\init(\bF^{(n)})$, is negligible for our purposes in this argument.
 Now by (\ref{eta-P2}) there are a bounded (in expectation) number of functions in $\rmod_{Y_0}(\bF^{(n)})$; 
 each function behaves essentially as the rescaled renewal process  associated with i.i.d.($\mu$)  summands, until 
 the wave ends, which happens within $O(1)$ steps when the function approaches the previous function.
We can apply the classical   invariance principle for renewal processes 
 (\cite{bill} Theorem 17.3) to show that
 (as $n \to \infty$) these functions behave as Brownian motion

This argument is sufficient to imply
\begin{equation}
\rmod_{Y_0}(\bF^{(n)}) \cd \rmod_{Y_0}(\bF^{(\infty)}) 
. \label{FnC3}
\end{equation}

\bigskip \noindent {\bf Step 4.}
To deduce  (\ref{FnC}) from  (\ref{FnC3}) we need to show
\begin{equation}
\rmod_{Y_0}(\bF^{(\infty)})  \cd \bF^{(\infty)} \mbox{ as } Y_0 \to \infty 
\label{Ytight}
\end{equation}
and the following analogous result for the sequence $(\bF^{(n)})$.  
Take $0 < Y_1 < Y_0 < \infty$ and $T_0 > 0$. 
Define an event $A(n,Y_1,Y_0,T_0)$ as
\begin{quote}
the random family $\rmod_{Y_0}(\bF^{(n)})$ and the random family $\bF^{(n)}$ do not coincide as regards 
functions $f$ with $f(0) \in [-Y_1,Y_1]$ considered only on $0 \le t \le T_0$.
\end{quote}
We then need to show: for each $Y_1, T_0$
\begin{equation}
\lim_{Y_0 \to \infty} \limsup_{n \to \infty} \Pr(A(n,Y_1,Y_0,T_0)) 
= 0 . 
\label{Ytight2}
\end{equation}
It is easy to prove  (\ref{Ytight})  directly from  basic properties of standard CBM.
Fix $Y_1, T_0$.  As in Figure 4, the particle positions of standard CBM at time $T_0$ determine
a point process $(t_u, u \in \Ints)$ on $\Reals$, each time-$T_0$ particle being the 
coalescence of all the  initially in the interval  $(t_{u-1},t_u)$.  So there are some random 
$ - \infty < u_L < u_R < \infty$ such that 
$ [-Y_1,Y_1] \subset [t_{u_L}, t_{u_R}]$, 
and then on the event 
$ [t_{u_L}, t_{u_R}] \subset [-Y_0, Y_0]$ 
we have, using the natural coupling of CBM processes, 
that the behavior of the particles of CBM$(\tilde{\eta}^{0},0)$ started within $ [-Y_1,Y_1]$
is unaffected by the behavior of the particles started outside $ [-Y_0,Y_0]$.
This is enough to prove  (\ref{Ytight}).

Unfortunately the simple argument above does not work for the analogous result (\ref{Ytight2}).
But there is a more crude argument which does work; for notational simplicity we will say this argument in the context 
of CBM but it then readily extends to the family  $\bF^{(n)}$.
We can construct  CBM started with particles at positions $\tilde{\eta}^{0}$ by first taking independent Brownian 
motions over $0 \le t < \infty$ started at each point of $\tilde{\eta}^{0}$ (call these virtual paths) and then, when two particles meet, we implement the coalescent by saying that the future path follows the virtual path of the particle with lower starting position. 
See the top panel of Figure 7, showing (lighter lines) the virtual paths of particles 
$a,b,\ldots,j$ and (darker lines) the coalescing process.
In words, we want to show the following.
\begin{quote}
(*) Fix $T_0$.  Consider CBM started with particles at positions 
$\tilde{\eta}^{0} \cap  [0,\infty]$ only. 
Then, in the natural coupling (induced by the construction above)  with CBM started with particles at all  positions 
$\tilde{\eta}^{0}$, the two process coincide, over $0 \le t \le T_0$, for the particles started on some random interval 
$[Z, \infty)$.
\end{quote}

The argument is illustrated in Figure 7, where particle $d$ is the first particle whose initial position is in 
$[0,\infty)$.
The top panel shows CBM started with particles at all  positions 
$\tilde{\eta}^{0}$, and the middle panel shows ``the modification"  of CBM started with particles at positions 
$\tilde{\eta}^{0} \cap  [0,\infty]$ only. 
The ``modification" affects (over the time interval shown)  particles $d,e,f,g,h$ but not subsequent particles.
We can upper bound the range of ``affects" as follows.
Consider only the ``virtual  paths", as shown in the top panel, and color path segments according to the
following rule, illustrated in the bottom panel.
Start to color the virtual path of particle $d$; when that path meets another virtual path, start coloring 
that path too; whenever a colored path meets an uncolored path, start coloring that path also.
This produces a branching network, and at time $T_0$ the paths associated with some set of particles have been colored.
It is easy to check that any particle not in that colored set is unaffected by the modification.
One can now check that only a finite mean number of particles will be colored 
(note this will include some particles started below $0$, for instance $b, c$ in the figure) and then that 
(*) holds with $\Ex Z < \infty$.

The point is that the calculation of a bound on $\Ex Z$ carries over to the rescaled family $(\bF^{(n)})$ because the 
``virtual paths" are just the centered renewal processes, which behave like Brownian motion; and we can use 
(\ref{Fnab2}) to show that only a bounded (in expectation)  number of functions enter the picture.

\newpage

\vspace*{-2.3in}

\includegraphics[width=4.0in]{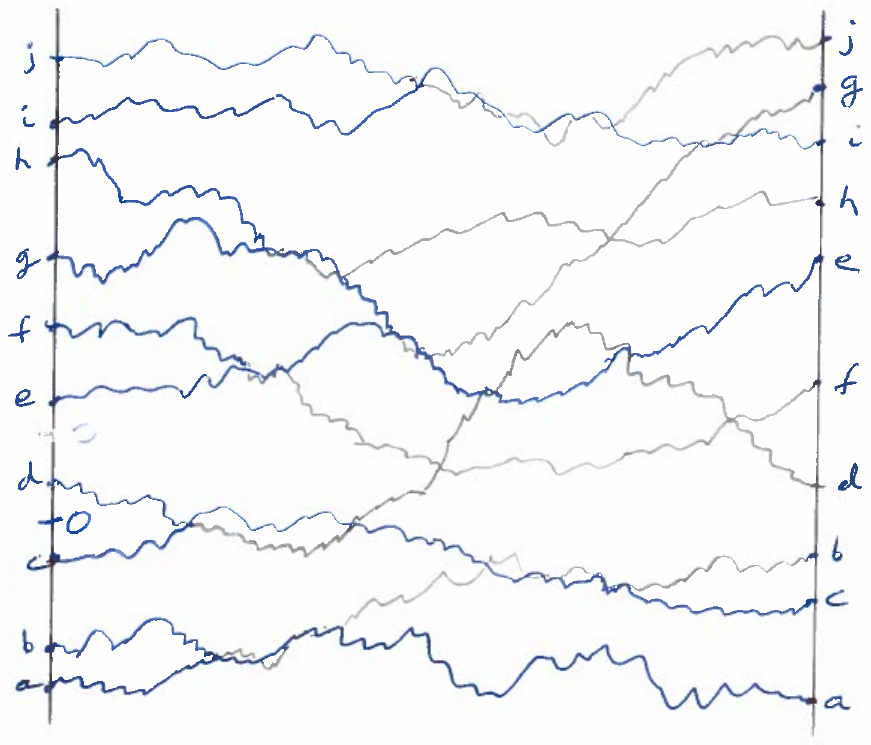}

\vspace*{-2.0in}

\includegraphics[width=4.0in]{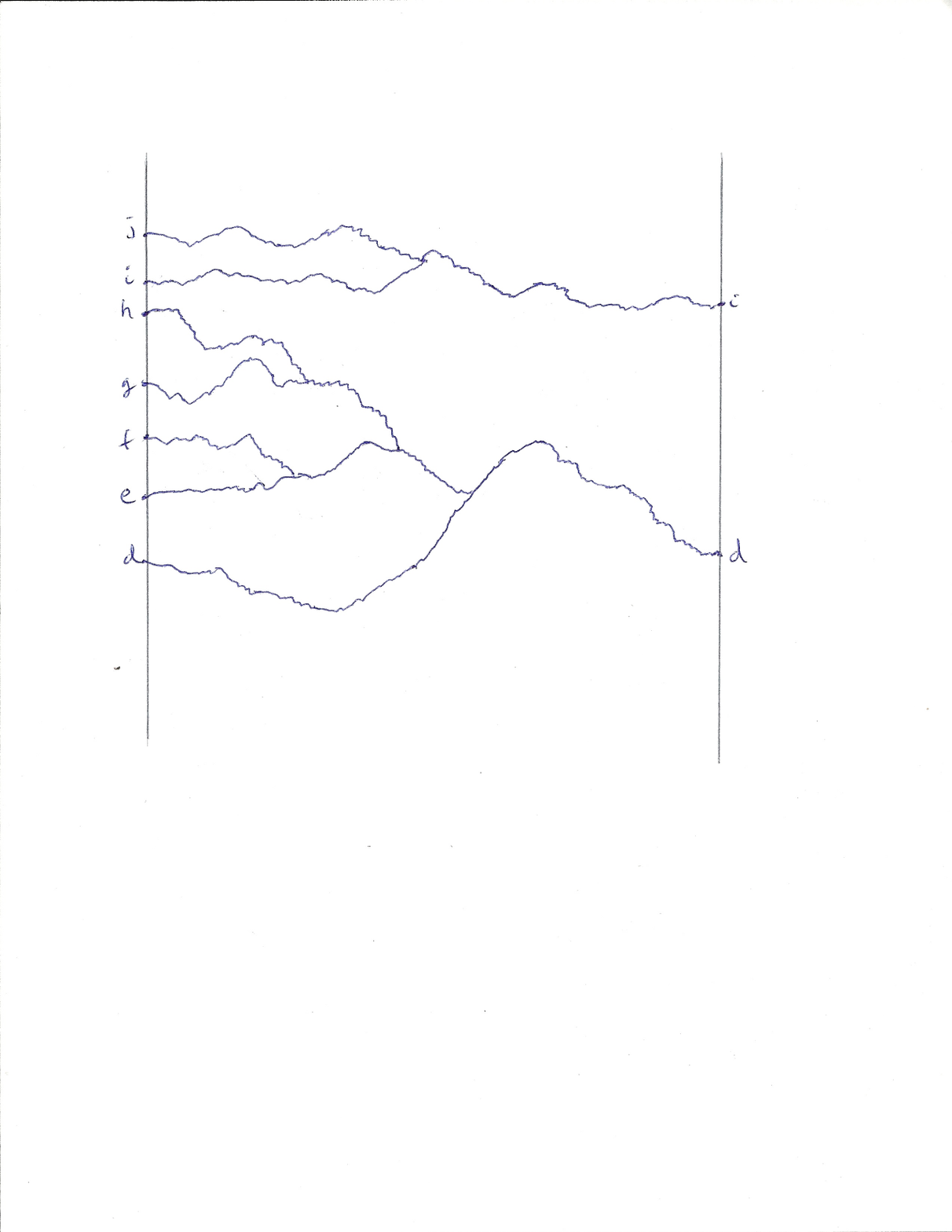}

\vspace*{-2.4in}

\includegraphics[width=4.0in]{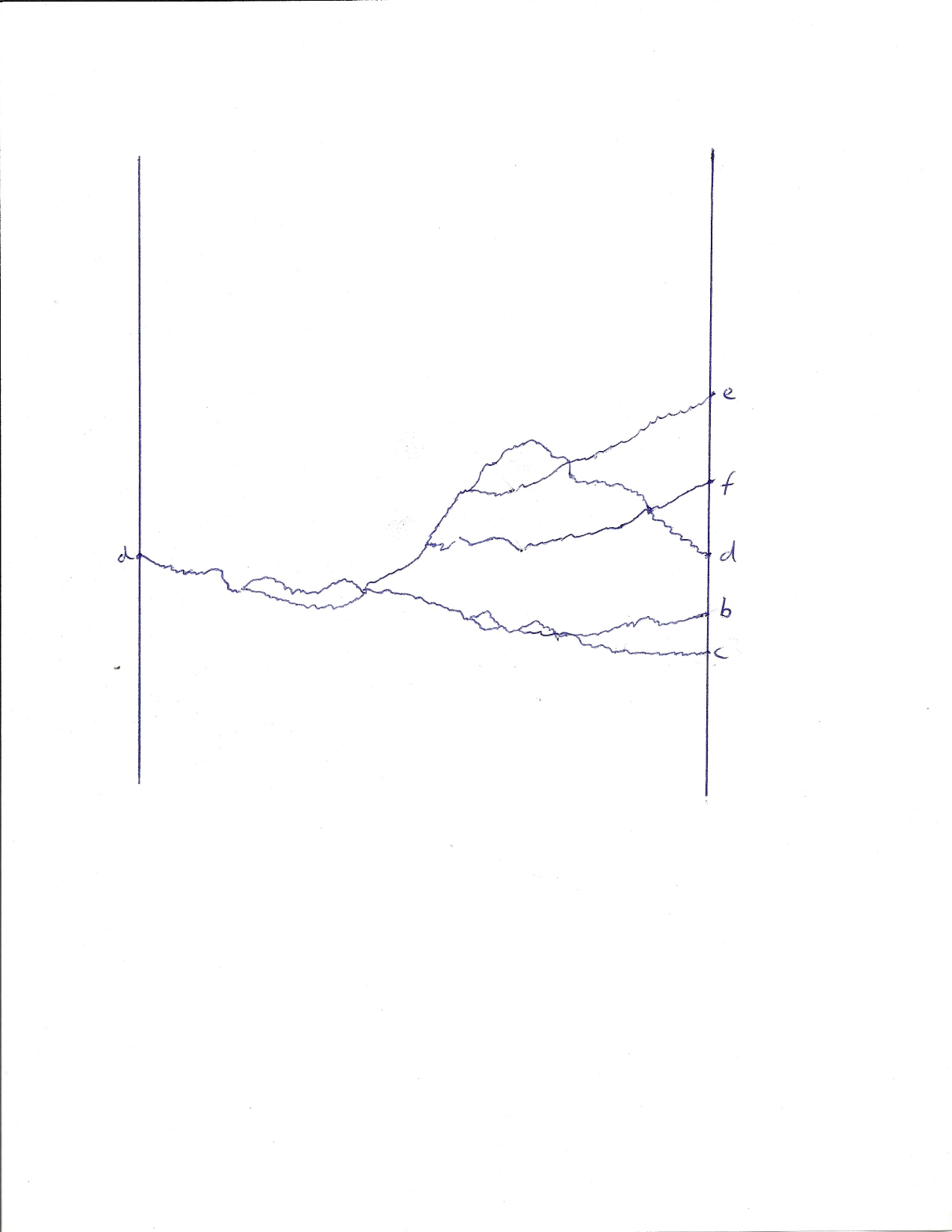}

\vspace*{-1.7in}

{\bf Figure 7.}  The construction in Step 4.

\newpage

\subsection{Consequences for wave lengths and space-time trajectories}
\label{sec:s-t}
Here we show what Theorem \ref{T:main} implies about the queue process,  outlining proofs of 
Theorem \ref{T:main_rate} and Corollary \ref{C:2C} stated  in section \ref{sec:what}.

\paragraph{Wave lengths.}
First let us switch from measuring wave length by rank to measuring it by position.  
That is, the length of the wave that creates $\bX(t)$ from $\bX(t-1)$ was defined implicitly at (\ref{Xit3}) as the rank
\[ W(t): = \min \{i: X_i(t) = X_{i-1}(t-1) \} \] 
of the first unmoved customer, 
whereas now we want to consider it as the position of that customer:
\[ L(t) := X_{W(t)}(t) . \]
Directly from the section \ref{sec:graphics} discussion we see that the length of wave at time $t$ is the position $x$ 
where $G(t,\cdot)$ coalesces with $G(t-1,\cdot)$,
that is where they first make the same jump.
Define a point process $\zeta^{(n)}$ to be the set of points
\[ \left\{ \sfrac{  t - U_{\tau_n}}{\sigma n^{1/2}}   \  : \ L(t) > n \right\}  .\]
Then, if $\tau_n \to \infty$ sufficiently fast, we have from Theorem \ref{T:main}
\begin{Corollary}
As $n \to \infty$ the point processes $\zeta^{(n)}$ converge in distribution to the spatial point process
$(B(y,1), y \in \Reals)$ of time-$1$ positions of particles in the standard CBM process.
\end{Corollary}
Now Lemma \ref{L:XkU} tells us that at a typical time in the queue process, the rank-$n$ customer is at position 
$n \pm O(n^{1/2})$.
Combining that result with the fact that, between waves of that length, the customer's position does not change and their rank 
decreases only by $O(n^{1/2})$, it is not hard to deduce the same result for wave lengths measured by rank:
if $\tau_n \to \infty$ sufficiently fast, then
\begin{Corollary}
\label{C:rateU}
As $n \to \infty$ the point processes 
$ \left\{ \sfrac{  t - U_{\tau_n}}{\sigma n^{1/2}}   \  : \ W(t) > n \right\}  $
 converge in distribution to the spatial point process
$\{B(y,1), y \in \Reals\}$ of time-$1$ positions of particles in the standard CBM process.
\end{Corollary}
This easily implies the ``rate" result stated as Theorem \ref{T:main_rate}.

\paragraph{Space-time trajectories.}
As defined, $\bF^{(n)}$ consists of the functions 
\[ t \to \widetilde{H}^{(n)}(y,t),  \quad 1 \le t < \infty  \]
as $y$ varies.
The same such function $f \in \bF^{(n)}$ arises from all $y$ in some interval, say the interval $[y_-(f),y_+(f)]$. 
Similarly, a function $f \in \CB$ is the function $t \to B(y,t), \ t \ge 1$ for all $y$ in some interval, which again we may write as the interval $[y_-(f),y_+(f)]$. 
Now note that the definition (\ref{def-Ht}) of  $\widetilde{H}^{(n)}(y,t) $ makes sense for $0 \le t \le 1$ also.  
For any fixed small $\delta > 0$, Theorem \ref{T:main} extends to showing convergence in distribution of 
$(t \to \widetilde{H}^{(n)}(y,t),  \ \delta \le t < \infty )$ 
to standard CBM over the interval $[\delta,\infty)$. 
This procedure enables us, by taking $\delta \downarrow 0$ ultimately,
to show that the Theorem \ref{T:main} convergence $\bF^{(n)} \cd   \CB$ extends to 
the ``augmented" setting  where the functions are marked by the intervals $[y_-(f),y_+(f)]$.
(At the technical level, this argument allows us to work within the space $\Fspace$ of locally finite collections of functions, instead
of some more complicated space needed to handle standard CBM on $(0,\infty)$).

We can now add some details to the verbal argument in section \ref{sec:what} for Corollary \ref{C:2C}.
Regarding the first component,
 the central point is as follows.
Take lines $\ell_n$ through $(1,0)$ of slope $- n^{1/2}$.
Write $(f^{(n)}_u, \ u \in \Ints)$ for the distinct functions in $\bF^{(n)}$  intercepted by line $\ell_n$, at some points 
$(x_u,s_u)$, and write $(y_-(f^{(n)}_u), y_+(f^{(n)}_u))$ for the interval of initial values of functions which coalesce with $f^{(n)}_u$ before $x_u$.
These represent the normalized move times of an individual chosen as being at position $n$ at time $U_{\tau_n}$
We know $\bF^{(n)} \cd   \CB$ in the ``augmented" sense above. 
This implies that the
associated interval endpoints  $(y_+(f^{(n)}_u) ,  u \in \Ints)$ (that is, the  normalized move times above)
converge in distribution as point processes 
to $(w_u, \ u \in \Ints)$, the time-$1$ positions of particles in CBM.
That is the assertion of Corollary \ref{C:2C}, as regards the first component.

Regarding the second component, the verbal argument in section \ref{sec:what} can be formalized to
show that for the process $(w_u, \ u \in \Ints)$ of positions of time-1 particles in CBM, the 
inter-particle distances are the rescaled limit of the numbers of customers who cross the starting position of the individual under consideration when that individual moves.
To then derive the stated result, the central issue is to identify ``numbers of customers who cross" with 
``distance the distinguished individual moves", which will follow from the lemma below.
Consider a typical wave which extends past position $n$ at some large time $t$.  
For each position $x$ in the spatial interval  $[n \pm B n^{1/2} ]$ there is some number 
 $M^{(n)}(x)$ of customers who cross position $x$ in the wave, and each individual $\iota$ in the spatial interval moves some distance $D^{(n)}(\iota)$.   From the CBM limit we know  
 these quantities are order $n^{1/2}$. 
 It is enough to show that all these customers move the same distance (to first order), in the following sense.
\begin{Lemma}
There exist random variables $M^{(n)}$ such that, for each fixed $B$,
  \[  \max_x  |  n^{-1/2}M^{(n)}(x) - M^{(n)} |  \to_p 0, \quad  \max_\iota  | n^{-1/2} D^{(n)}(\iota) - M^{(n)} |  \to_p 0  \mbox{ as } n \to \infty \]
 where the maxima are over  the spatial interval  $[n \pm B n^{1/2} ]$.
 \end{Lemma}
Consider the
 inter-customer distances $(\xi_i(t-1), \xi_i(t), i \ge 1)$ immediately before and after the wave at time $t$.
If, over the interval $[n \pm B n^{1/2} ]$ under consideration, these were i.i.d.($\mu$) then the assertion of the lemma would be clear, 
using (for instance)  large deviation bounds for i.i.d. sums.  
Analyzing wave lengths in sections \ref{sec:UB} and \ref{sec:LB} was complicated because we did not know orders of magnitude.
Now we know the power law $\Pr(W > w) \sim c w^{-1/2}$ for wave length, so  conditional on a wave entering the interval 
$[n \pm B n^{1/2} ]$ the probability it goes only $o(n)$ further is $o(1)$.  
From the CBM limit we know that the first customer in the interval moves a distance $M^{(n)}$ of order $n^{1/2}$.
Given that distance, the subsequent inter-customer distances $\xi_i(t)$ over the interval are  i.i.d.($\mu$) 
conditioned on a event of probability $1 - o(1)$.
The same was true for the previous wave that created the inter-customer distances $\xi_i(t-1)$.
So the conclusion of the lemma follows from the i.i.d. case.

\newpage
\section{Final remarks}
\subsection{Some data}
\label{sec:OAK}
Spending 17 minutes in line at security at Oakland airport enabled collection of the data in Figure 8, 
whose points  show rank and time after each move. 
The  approximate straight line reflects the fact that people are being served at an approximately constant rate.
The initially larger gaps between points indicate the ``wave" phenomenon studied in this paper.

\setlength{\unitlength}{0.0024in}
\begin{picture}(1850,1150)
\put(0,0){\line(1,0){1800}}
\put(0,0){\line(0,1){1060}}
\put(1760,0){\circle*{15}}
\put(1690,28){\circle*{15}}
\put(1640,69){\circle*{15}}
\put(1560,102){\circle*{15}}
\put(1480,180){\circle*{15}}
\put(1450,211){\circle*{15}}
\put(1420,238){\circle*{15}}
\put(1370,279){\circle*{15}}
\put(1310,330){\circle*{15}}
\put(1260,354){\circle*{15}}
\put(1220,405){\circle*{15}}
\put(1200,441){\circle*{15}}
\put(1140,451){\circle*{15}}
\put(1110,481){\circle*{15}}
\put(1070,490){\circle*{15}}
\put(1040,512){\circle*{15}}
\put(1000,526){\circle*{15}}
\put(920,543){\circle*{15}}
\put(870,562){\circle*{15}}
\put(820,588){\circle*{15}}
\put(790,614){\circle*{15}}
\put(770,631){\circle*{15}}
\put(710,641){\circle*{15}}
\put(660,650){\circle*{15}}
\put(620,684){\circle*{15}}
\put(560,700){\circle*{15}}
\put(540,726){\circle*{15}}
\put(510,741){\circle*{15}}
\put(490,758){\circle*{15}}
\put(450,776){\circle*{15}}
\put(390,810){\circle*{15}}
\put(370,833){\circle*{15}}
\put(340,844){\circle*{15}}
\put(330,855){\circle*{15}}
\put(290,861){\circle*{15}}
\put(270,890){\circle*{15}}
\put(230,912){\circle*{15}}
\put(200,932){\circle*{15}}
\put(190,948){\circle*{15}}
\put(170,963){\circle*{15}}
\put(160,973){\circle*{15}}
\put(130,984){\circle*{15}}
\put(100,1003){\circle*{15}}
\put(60,1014){\circle*{15}}
\put(50,1026){\circle*{15}}
\put(30,1034){\circle*{15}}
\multiput(200,0)(200,0){8}{\line(0,-1){23}}
\multiput(0,120)(0,120){8}{\line(-1,0){22}}
\put(170,-77){20}
\put(370,-77){40}
\put(570,-77){60}
\put(770,-77){80}
\put(960,-77){100}
\put(1160,-77){120}
\put(1360,-77){140}
\put(1560,-77){160}
\put(-70,100){2}
\put(-70,220){4}
\put(-70,340){6}
\put(-70,460){8}
\put(-90,580){10}
\put(-90,700){12}
\put(-90,820){14}
\put(-90,940){16}
\put(740,-140){rank in line}
\put(-250,530){time}
\put(-270, 460){(min.)}
\end{picture}

\vspace*{0.4in}

{\bf Figure 8.}  Data: progress in a long queue.

\subsection{Stationarity}

As mentioned earlier, 
it is possible that continuing the argument for Lemma \ref{L:almost} would give a ``coupling from the past" proof of 
convergence to a stationary distribution (Conjecture \ref{C-1}).
Assuming existence of a stationary distribution,  a stronger version  of Lemma \ref{L:almost} would be
\begin{OP}
Suppose $\mu$ has a density that is bounded below on $[c_-,c^+]$ and suppose 
the  stationary distribution $\bX(\infty)$ exists.  
Is it true that
the joint distribution of all inter-customer distances 
$(X_i(\infty) - X_{i-1}(\infty), i \ge 1)$
 is absolutely continuous with respect to infinite product measure $\mu^\infty$.
\end{OP}
Heuristics based on the the Kakutani equivalence theorem suggest the answer is ``no".

Other questions involve mixing or ``burn in": how long does it take, from an arbitrary start, for the distribution of the 
process restricted to the first $k$ positions be close to the stationary distribution?
From the initial configuration with all inter-customer distances equal to $c_-$, 
the time must be at least  order $k$: is this the correct worst-case-start order?

\subsection{Robustness to model generalizations}
\label{sec:robust}
Let us amplify an earlier statement   ``intuitively the model behavior should also be robust to different customers having different parameters 
$\mu, c_-,c^+$". 
What we have  in mind is a prior distribution over these parameters, and the initial time-0 queue of customers 
having parameters sampled i.i.d. from the prior.
The key issue is that the analog of the symmetrized random walk at (\ref{Ssym}),
representing the same block of customers at different times,
 is still 
(even conditioned on the realization of parameters) a sum of independent mean-zero increments having 
 (under mild assumptions) the same asymptotic behavior as given in Lemma \ref{LRW}.  
The proof should carry over without essential changes.

\subsection{Finite queues}
A natural alternative to our infinite-queue setup would be to take Poisson arrivals with (subcritical) rate 
$\rho < 1$.  Because such a $M/D/1$ queue is empty for a non-vanishing proportion of time, a stationary distribution clearly exists.  
However, to state the analog of Theorem \ref{T:main} for the stationary process one still needs a double limit 
(rank $j \to \infty$ and traffic intensity $\rho \uparrow 1$); 
the infinite-queue setup should, if one can prove Conjecture \ref{C-1}, allow  a single limit statement of form 
(\ref{single-limit}).

\subsection{Other possible proof techniques}
As mentioned earlier we suspect there may be simpler proofs using other techniques.
Here are brief comments on possible  proof techniques we did think about but did not use.

\paragraph{Forwards coupling.}
The classical notion of (forwards) Markov chain coupling -- showing that versions of the process from two different initial states 
can be coupled so that they eventually coincide in some useful sense -- does not seem to work for our process.

\paragraph{Feller property.}
A Feller process $(X(t))$ on a compact space always has at least one stationary distribution, 
by considering subsequential weak limits of $X(U_\tau)$. 
Our spatial queue process does have a compact state space $\mathbb{X}$; it is not precisely Feller, because of 
the discontinuity in the construction (\ref{Xit1} -  \ref{Xit3}),
but  Feller continuity holds a.e. with respect to product uniform measure, so under mild conditions 
(e.g. that $\mu$ has a density component) the same argument applies to prove existence of at least one stationary distribution.
However, such ``soft" arguments cannot establish uniqueness of stationary distribution  or convergence thereto. 
One could rephrase our limit results in terms of stationary processes with such subsequential stationary distributions,
 but that hardly seems a simpler reformulation.

\paragraph{CBM as a space-indexed Markov process.}
A conceptually different starting point for an  analysis of our spatial queue model is to observe that the  standard CBM process 
$(B_y(s), 0 \le s < \infty, y \in \Reals)$ with $B_y(0) = y$ 
of section \ref{sec:CBM}, with ``time $s$" and ``space" $y$, 
can in fact be viewed in the opposite way.  Define
\[ X_y = (B_y(s) - y, \ 0 \le s < \infty) \]
so that $X_y$ takes values in the space $C_0(\Reals^+)$ of continuous functions $f$ with $f(0)=0$. 
Now the process $(X_y, - \infty < y < \infty)$ with ``time" $y$ is a continuous-time $C_0(\Reals^+)$-valued Markov process.
Apparently CBM has not been studied explicitly in this way.  
In principle one could determine its generator and seek to apply general techniques for weak convergence 
of discrete-time Markov chains to continuous-time limits.  
But we have not attempted this approach.

\subsection{Generalizations of CBM}
Finally we mention that the  extensions  of CBM to the {\em Brownian web} and {\em Brownian net} have been  shown to arise as scaling limits of various one-dimensional models \cite{Bweb}. 
But our queue model seems to  be a novel addition to this collection of CBM-related models.


\paragraph{Acknowledgments.}
I thank Weijian Han and Dan Lanoue for  simulations of the process.

\newpage

\begin{thebibliography}{10}

\bibitem{me-RW}
D.~Aldous.
\newblock Overview of probability in the real world project, 2015.
\newblock web page http://www.stat.berkeley.edu/~aldous/Real-World/cover.html.

\bibitem{arratia1979}
Richard~Alejandro Arratia.
\newblock {\em Coalescing Brownian motions on the line}.
\newblock PhD thesis, University of Wisconsin--Madison, 1979.

\bibitem{bill}
Patrick Billingsley.
\newblock {\em Convergence of probability measures}.
\newblock John Wiley \& Sons, Inc., New York-London-Sydney, 1968.

\bibitem{blank_2010}
Michael Blank.
\newblock Metric properties of discrete time exclusion type processes in
  continuum.
\newblock {\em J. Stat. Phys.}, 140(1):170--197, 2010.

\bibitem{borodin}
Alexei Borodin, Patrik~L. Ferrari, Michael Pr{\"a}hofer, and Tomohiro Sasamoto.
\newblock Fluctuation properties of the {TASEP} with periodic initial
  configuration.
\newblock {\em J. Stat. Phys.}, 129(5-6):1055--1080, 2007.

\bibitem{doi:10.3141/2177-08}
Darcy Bullock, Ross Haseman, Jason Wasson, and Robert Spitler.
\newblock Automated measurement of wait times at airport security.
\newblock {\em Transportation Research Record: Journal of the Transportation
  Research Board}, 2177:60--68, 2010.

\bibitem{gray}
Lawrence Gray and David Griffeath.
\newblock The ergodic theory of traffic jams.
\newblock {\em J. Statist. Phys.}, 105(3-4):413--452, 2001.

\bibitem{kleinrock1975queueing}
Leonard Kleinrock.
\newblock {\em Queueing Systems. Volume 1: Theory}.
\newblock Wiley-Interscience, 1975.

\bibitem{Bweb}
E.~{Schertzer}, R.~{Sun}, and J.~M. {Swart}.
\newblock {The Brownian web, the Brownian net, and their universality}.
\newblock {\em ArXiv 1506.00724}, June 2015.

\bibitem{toth-werner}
B{\'a}lint T{\'o}th and Wendelin Werner.
\newblock The true self-repelling motion.
\newblock {\em Probab. Theory Related Fields}, 111(3):375--452, 1998.

\end{thebibliography}

\end{document}